\RequirePackage{fix-cm}
\documentclass{article}
\usepackage[utf8]{inputenc} 
\usepackage{amssymb,amsmath} 
\usepackage{graphicx} 
\usepackage{eurosym} 
\usepackage{xcolor} 
\usepackage{booktabs}
\usepackage{pifont} 
\usepackage{algorithm} 
\usepackage{algpseudocode}
\usepackage{multicol}
\usepackage{multirow}
\usepackage{subcaption}
\usepackage{hyperref}
\usepackage{setspace}
\usepackage[title]{appendix}
\usepackage[square,numbers]{natbib}
\usepackage[affil-it]{authblk}
\usepackage{url}

\begin{document}
	
	\title{A two-phase stochastic programming approach to biomass supply planning for combined heat and power plants}

\author[1]{Ignacio~Blanco}
\author[1,*]{Daniela Guericke}
\author[2]{Juan M. Morales}
\author[1]{Henrik~Madsen}

\affil[1]{\small Technical University of Denmark, Department for Applied Mathematics and Computer Science, Richard Petersens Plads, 2800 Kgs. Lyngby, Denmark}
\affil[2]{\small M\'alaga University, Department of Applied Mathematics, Escuela de Ingenierías Industriales, 29071 M\'alaga, Spain}
\affil[*]{\small Corresponding author: Daniela Guericke, dngk@dtu.dk}

	\maketitle
		
	\begin{abstract}
			 Due to the new carbon neutral policies, many district heating operators start operating their combined heat and power (CHP) plants using different types of biomass instead of fossil fuel. The contracts with the biomass suppliers are negotiated months in advance and involve many uncertainties from the energy producer’s side. The demand for biomass is uncertain at that time, and heat demand and electricity prices vary drastically during the planning period. Furthermore, the optimal operation of combined heat and power plants has to consider the existing synergies between the power and heating systems.
		We propose a solution method using stochastic optimization to support the biomass supply planning for combined heat and power plants. Our two-phase approach determines mid-term decisions about biomass supply contracts as well as short-term decisions regarding the optimal production of the producer to ensure profitability and feasibility. We present  results based on two realistic test cases.
		
	\end{abstract}

\section{Introduction}
The integration of different energy systems is one step towards a fossil-free energy system, which many developed countries target today. By integrating different energy systems, such us heat and power, a higher share of volatile renewable energies, e.g., wind energy, can be used efficiently \cite{LUND2011419}.
In areas with large district heating networks, one way to achieve this integration is using combined heat and power (CHP) plants that produce heat and power simultaneously. By co-optimizing the production of both, the efficiency of the system is increased while providing flexibility to the power grid and satisfying the heat demand in the district heating network. 
Due to the neutral carbon policies imposed by the authorities, a shift from traditional fuels to renewable resources is taking place. Denmark has a widespread use of district heating and CHP plants and the government supports the use of biomass to produce heat and power. With subsidies and tax benefits, it has become profitable for large-scale CHP plants to change from, e.g., coal or natural gas to biomass \cite{197363:online}.

The use of biomass as fuel for CHP plants raises some challenges in the planning of the supply and in the operation of the plant. Many different types of biomass are used to produce heat and power \cite{yin2008grate} but the most common type of biomass used for large-scale CHP producers is wood pellets. Due to their high energy content, wood pellets facilitate a more efficient transport because smaller volumes are required. In addition, the low moisture content of wood pellets allows a better conservation of the product resulting in a larger storage capacity \cite{senechal2009logistic}. In combination with neutral carbon policy incentives for biomass, the wood pellet is becoming a candidate to substitute coal in CHP plants. However, comparing the supply of wood pellets, or biomass in general, with supply of natural gas, the former has some disadvantages. First, natural gas prices have been dropping since 2008 and, second, it exists a well-developed infrastructure for natural gas, which allows the producer to be directly connected to the gas network. On the contrary, biomass is transported long ways and contracts with the supplier must be agreed beforehand for a long horizon (one to three years) involving a high degree of uncertainty at negotiation because the final amount  is  unknown. It is crucial for CHP operators to optimize their biomass contracts to be competitive with gas-fired plants. 

In this work, we propose a solution approach based on stochastic programming \cite{birge2011introduction} to optimize the yearly biomass contracting decisions for a CHP operator taking into account the uncertainty at the point of negotiation.  Furthermore, the approach also determines the optimal operation of the plant to maximize profits and satisfy the heat demand on weekly basis throughout the year. 
\section{Literature review}
Several models for the optimal operation of CHP systems, where different aspects of the problem are  highlighted, have been proposed. We refer for example to \cite{aringhieri2003optimal,rolfsman2004combined,dvovrak2012combined,christidis2012contribution,Rong2007,pirouti2011optimal}. These solution approaches determine the optimal production of both commodities (heat and power) at different levels of detail, but do not consider uncertainties and supply contracts for fuel explicitly.

Since then several approaches that apply stochastic programming for the operational planning were developed.
\cite{alipour2014stochastic} solve the operational scheduling for an industrial customer that owns an integrated system formed by CHP units, conventional power production and heat only units. The method uses electricity market sales and demand response programs to integrate the uncertainty caused by electricity prices and load. 
An optimal operation of a portfolio of different CHP systems in a district heating network is studied in \cite{nielsen2016economic}. The authors consider uncertain heat demand and electricity prices and show that the system profits from leveraging a thermal storage to handle this uncertainty.  \cite{Kumbartzky2017} present a multi-stage stochastic program for optimizing the operation of a gas-fired CHP plant and deriving bids for the German spot and balancing markets. The considered uncertainty are electricity prices. 
\cite{dimoulkas2014constructing} propose a stochastic program including technical aspects of a extraction-condensing CHP plant for optimizing the hourly operation under price and demand uncertainty. The authors use this model to determine bidding curves for the day-ahead  market. In \cite{dimoulkas2015probabilistic} this model is revisited with more focus on the joint production scheduling of two CHP plants. The operational planning problem in our work is similar to these two formulations, but extended with further characteristics regarding the biomass contracts deliveries and technical constraints. 

The above mentioned publications assume instantaneous fuel supply and, therefore, do not consider fuel supply decisions. Another stream of publications explicitly concentrates on the biomass supply chain planning for power generation considering processing of biomass, transport and logistics aspects.
The OPTIMASS model for strategic and tactical biomass supply chain planning is presented in \cite{DeMeyer2015}. The formulation is based on a facility location planning problem that includes the processing of the biomass to determine locations and capacities of facilities in the supply chain and allocation of biomass sites to conversion facilities. The final usage of biomass in electricity production is not part of this study. \cite{frombo2009decision} present a decision support system for a forest biomass supply chain deciding on the locations and capacities as well as assignment of biomass sources to power plants. \cite{QUDDUS201827} present a two-stage stochastic program with chance constraints for biomass supply chain planning under biomass availability uncertainty. The demand is based on markets and not single plants. \cite{bruglieri2008optimal} model the biomass-based energy production process, which includes deciding the location of plants as well as flow and conversion of commodities where one commodity is electricity. The model focuses on long-term decisions.  

In this work the perspective of a power plant that receives biomass from third party suppliers is considered. Furthermore, we investigate the integration of long-term biomass supply decisions with the operational planning of the production. Similar settings have been studied in the following publications.
\cite{maurovich2016optimal} consider the fuel supply of gas for a consumer having a micro CHP and a heat boiler. Their multi-stage stochastic program decides on how much gas to buy on the spot or the monthly and weekly futures market, while electricity can be sold with similar market instruments.  The model has a monthly planning horizon and abstracts from more detailed considerations regarding the operation of the system. 
In \cite{shabani2014tactical}, a general overview of the benefits of using stochastic programming to incorporate the uncertainty involved in the biomass supply chain for a power producer on a tactical planning level is given. The authors formulate a one year planning problem considering the amount of biomass supply from different suppliers, storage and the expected power production on a monthly basis. \cite{Jensen2017} consider the supply chain connected to a biogas CHP plant and use a network flow model formulation. The model includes conversion to biogas and production with a CHP or heat boiler as well as transportation costs.
\cite{SCOTT2015226} address biomass supplier selection combining an analytic hierarchy process (AHP) with a chance constraint program to address stakeholders and uncertainties in this setting. Their focus is ensuring the quality of the biomass by blending biomass from different kinds and suppliers to fulfill the overall demand. The solution approach disregards the production level and delivery times.
Finally, \cite{Chiarandini2014} use stochastic programming for optimal biomass contracting decisions in a long-term planning horizon. The model decides which biomass contracts should be settled with the suppliers. They model the contracts as well as the deliveries and production to provide a basis for this decision. Due to the planning horizon and short time periods, the model results in a computationally hard two-stage stochastic program. 

Our work differs from \cite{Chiarandini2014} regarding the modeling of contracting decisions and the overall solution approach. Delivery times and amounts for contracts in \cite{Chiarandini2014} are fixed and the decision-maker can just decide which contracts are selected. On the contrary, our approach allows more flexibility to decide on the amount to be supplied and the delivery time. As a consequence, the exact delivery time and precise quantity are determined once we are getting closer to the energy delivery. Furthermore, we reduce the computational complexity of the planning problem by presenting a two-phase approach. 

The main contributions of our work are the following:
\begin{enumerate}
	\item We propose a two-phase solution approach that combines biomass contracting decisions with the optimal operation of the CHP plant.  Therefore, it provides two models that can be used by an operator for long-term and operational planning, respectively. The first phase concentrates on the biomass contract selection at the beginning of the year considering production on a weekly less detailed basis and, therefore, reducing the complexity of the problem. The second phase optimizes the weekly operation of the system on a detailed hourly basis and takes the biomass contract decisions into account. The overall solution approach considers relevant technical requirements and resembles the planning process in practice.
	\item Our modeling of biomass contracts offers a high degree of flexibility. Completely fixed contracts can be investigated as well as more flexible contracts regarding  amounts of deliveries. We include the possibility to buy options on the biomass amount to be able to adjust the delivery quantity during the course of the year. This is a new feature whose benefits are worth of investigation, at least from the standpoint of a CHP producer. 
	\item Furthermore, we use a receding horizon approach to improve the results of our weekly operational planning, because it is important to take initial information from previous weeks into account and have a feasible transition. This also allows us to update the scenarios with new information.
\end{enumerate}
The remainder of this papers is organized as follows. A detailed description of the planning problem is given in Section \ref{ProblemDescription}.
Our solution approach and the respective model formulations are presented in Section \ref{Solutionapproach}. In Section \ref{casestudies}, we analyze two realistic case studies. The section includes a description of the data, experimental setup and scenario generation. The numerical results are stated in Section \ref{sec:results}. Finally, Section \ref{conclusion} summarizes our work and gives an outlook.  

\section{Problem description}
\label{ProblemDescription}
\begin{figure}
	\centering
	\includegraphics[width=0.7\columnwidth]{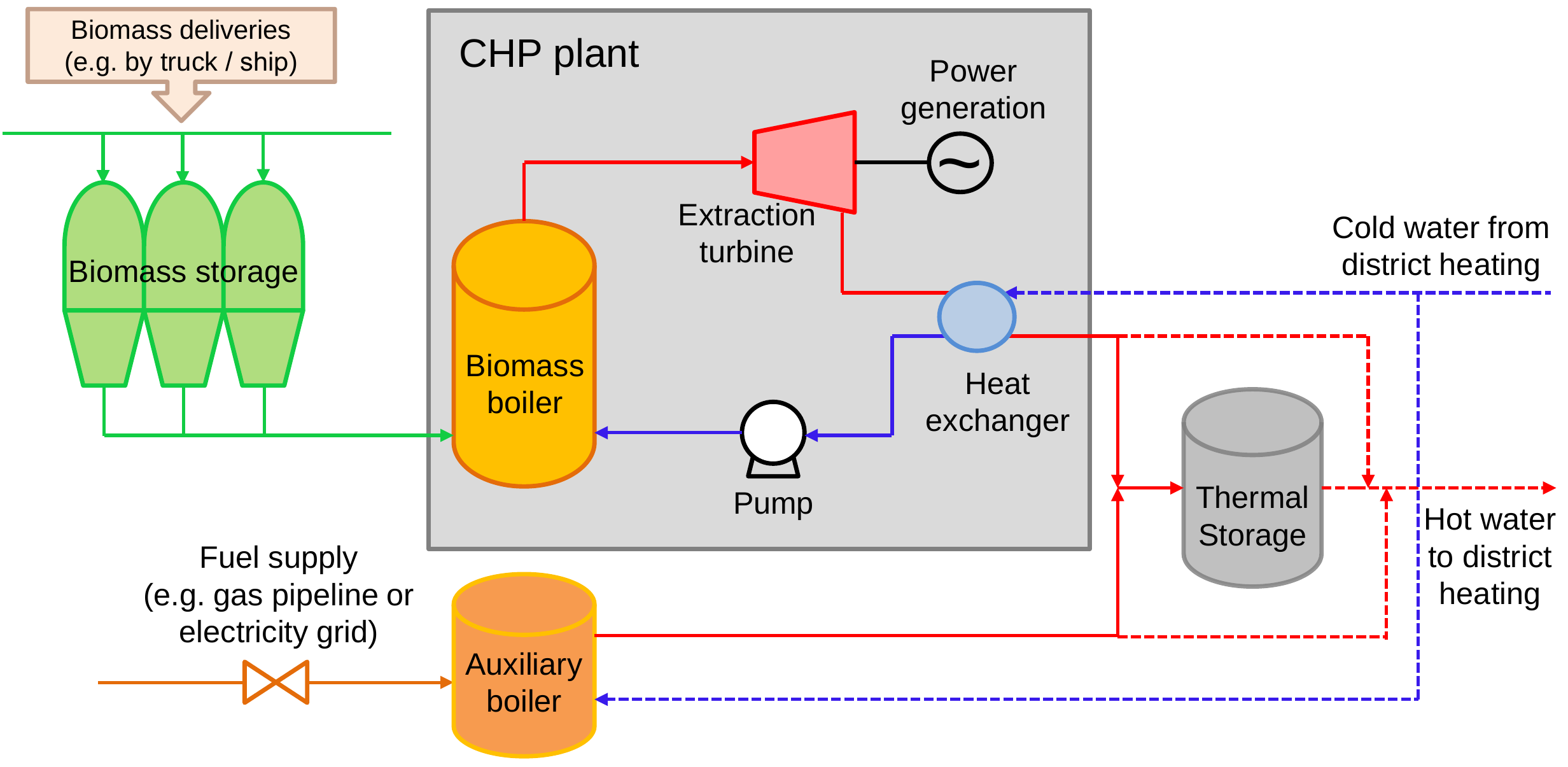}
	\caption{Overview of components in the planning problem}
	\label{fig:flowchart}
\end{figure}

In this section, we describe the biomass supply planning problem including used sets and parameters. For quick reference, we also provide an overview of parameters and sets in Table \ref{tab:parameters}.

An overview of the components in the planning problem   is given in Fig. \ref{fig:flowchart}.  We consider a power and heat producer directly connected to a district heating network. The producer operates a CHP plant fueled by biomass and an auxiliary heat producing unit (e.g. gas boiler, electric boiler or heat pump). Both units can supply the district heating network directly but are also connected to a thermal storage, which can store hot water for later heat supply.

The biomass delivered by suppliers according to the contracts is unloaded into the biomass storage and withdrawn from the storage for later use (i.e. no direct supply to the boiler). We assume that fuel for the additional heat-only unit is provided directly and instantaneously without storage and deliveries. This assumption stems from the setting of a gas boiler connected to the gas network or an electric boiler connected to the electricity grid.

In practice, biomass contracts are often agreed for a period of one year or more, defining the amount of biomass and a preliminary delivery schedule. The actual delivery time is revised in the course of the year. We model two different types of contracts, namely fixed and flexible. \textit{Fixed} contracts are cheaper but offer no possibility to alter the delivery amount afterward. \textit{Flexible} contracts are more expensive than fixed contracts, but the operator has the opportunity to buy an option of changing the  amount. In the beginning of the year, in addition to the delivery amount, the options for up- and/or down-scaling the amount are settled, while the producer has to pay extra for those options. The possibility of buying options to change the biomass delivery amount is a new concept that is studied in this paper. It provides the power producer with additional flexibility that can be beneficial especially in the long term when the actual demand is still uncertain. Also from the supplier's side this could be an interesting instrument, because it offers additional incomes from selling options while the amounts can be shifted between different customers. However, the supplier side is not the focus of this paper.

The input to our solution approach is a set of possible contracts $\mathcal{J}$, a set of scenarios $\Omega$ and a set of periods $\mathcal{T} = \lbrace 1, \ldots, |\mathcal{T}|\rbrace$. The first planning period is always denoted with 1, so that initial values are given values for period 0 (e.g. for storage levels $\delta_{0,\omega}$ and $s_{0,\omega}$). Each contract $j$ has a minimum and maximum amount per delivery ($\underline{B}_j,\overline{B}_j$), a minimum and maximum number of deliveries per planning horizon ($\underline{N}_j,\overline{N}_j$) and a minimum time between deliveries ($F_j$). If contract $j$ offers up-scaling and down-scaling options, the maximum limitations are given by $O^+_j$ and $O^-_j$ (in percent deviation from the nominal amount), respectively. For fixed contracts these parameters are set to zero ($O^+_j=O^-_j=0$). The cost for the fixed, up-scaling and down-scaling amount are given by $C^B_j, C^{B+}_j$ and $C^{B-}_j$, respectively. The cost are given per MWh, because the payment in practice is determined based on the energy content of the biomass in Gigajoule, which can be directly transformed to MWh. This means that the payment does not depend only on the amount in tonnes but also on the quality of the biomass, the so-called calorific value. Transportation costs are considered only indirectly, because the supplier has to cover these and can include them in the biomass cost per MWh. Furthermore, we assume that the supplier has the responsibility to deliver the contracted amount.
As mentioned above, the biomass is delivered to the biomass storage, which is limited by a minimum safety and maximum storage level ($\underline{\Delta}_t, \overline{\Delta}$). The initial storage level is given for period 0 and the outflow per period is restricted to a maximum of $\Delta^\text{F}$. To avoid congestion at the storage due to several deliveries at the same time, the time distance between deliveries must be at least $\Delta^\text{W}$ periods.

\begin{figure}
	\centering
	\includegraphics[width=0.35\columnwidth]{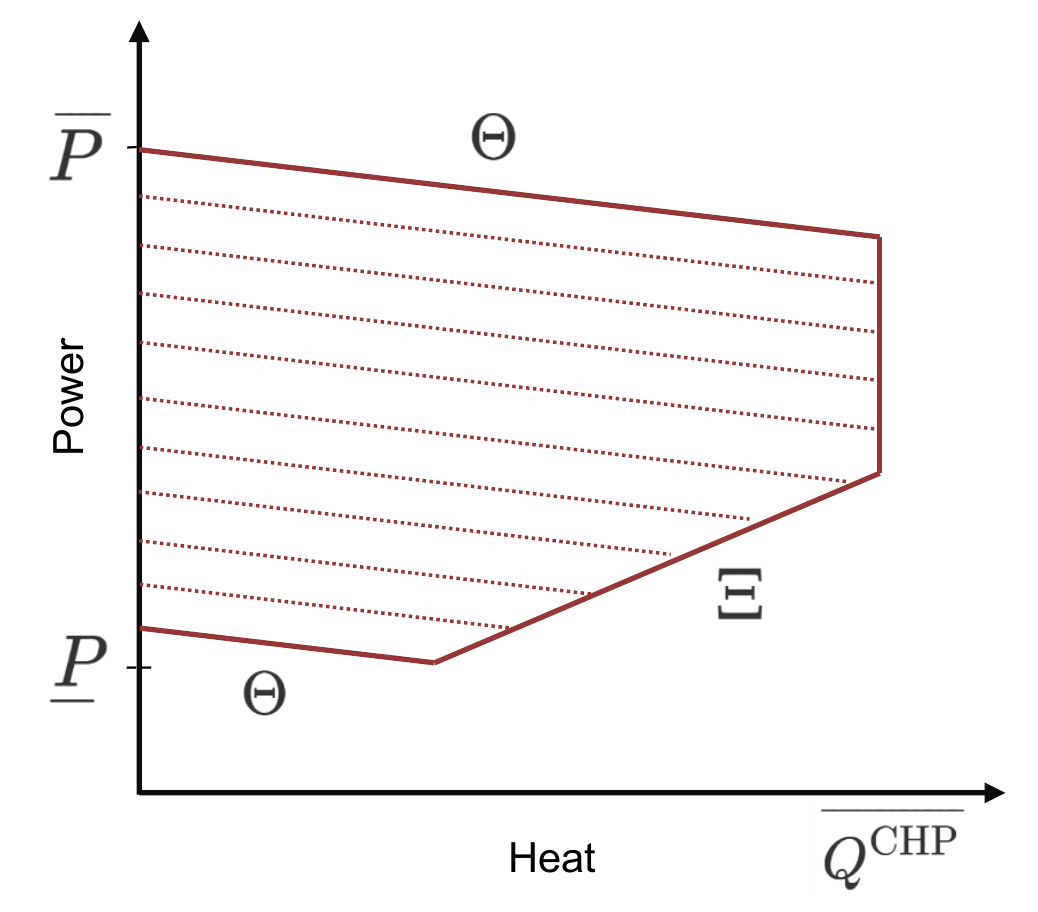}
	\caption{Feasible production region of extraction-condensing unit in CHP plant}
	\label{fig:chp_feasibleregion}
\end{figure}

\begin{table}
	\caption{Sets and parameters}
	\footnotesize
	\begin{tabular}{lp{0.85\columnwidth}}\toprule
		$\mathcal{J}$& Set of biomass contracts $j$\\
		$\mathcal{W}$& Set of weeks $\mathrm{w}$ \big($\mathcal{W} = \lbrace 1,\ldots, |\mathcal{W}| \rbrace$\big)\\
		$\mathcal{T}$& Set of time periods $t$\\
		$\mathcal{T}_\mathrm{w}$& Set of time periods $t$ in week $\mathrm{w}$\\
		$\Omega$& Set of scenarios $\omega$\\ \midrule
		$\pi_{\omega}$ &  Probability of scenario  $\omega$ \\
		$D_{t,\omega}$ & Heat demand in period $t$ in scenario $\omega$ [MWt/period] \\
		$L_{t,\omega}$ & Negative costs, i.e. profit, for selling electricity in period $t$ in scenario $\omega$ [\euro/MWe]\\
		$C^{\text{AUX}}_{t,\omega}$ & Operational cost of auxiliary boiler in period $t$ in scenario $\omega$ [\euro/MWt]\\
		$C^{\text{CHP}}$ & Operational cost of CHP plant [\euro/MWt]\\
		$C^{\text{SU}}   $ &        Start up cost for CHP   [\euro/MWt]\\
		$C^{\text{SD}} $ &           Shut down cost for CHP [\euro/MWt]\\
		$C^{\text{I}}$ & Inventory cost for biomass storage   [\euro/MWt]\\
		$C^{\text{F}}_{t,\omega}$ & Cost of fuel for auxiliary boiler in period $t$ and scenario $\omega$ [\euro/MWt]\\
		\vspace{0.1cm}
		$C^{\text{O\&M}}_{AUX}$ & Operational cost for auxiliary boiler [\euro/MWt]\\
		\vspace{0.1cm}
		$T^{\text{EP}}$ & Tax for electricity production [\euro/MWe]\\
		$T^{\text{AUX}}$ & Tax for production with auxiliary boiler [\euro/MWt]\\
		$T^{\text{CO}_2}$ & CO$_{2}$ emission tax [\euro/MWt]\\
		$C^{\text{B}}_{j}$ &  Cost for biomass in contract $j$      [\euro/MWt]\\
		$C^{\text{B}+}_{j}$ & Cost for up-scaling biomass amount in contract $j$    [\euro/MWt]\\
		$C^{\text{B}-}_{j}$ & Cost for down-scaling biomass amount in contract $j$   [\euro/MWt]\\
		$\underline{B_{j}}, \overline{B_{j}}$ &    Minimum/maximum amount biomass offered per delivery by contract $j$   [MWt]\\
		$\underline{N_j}, \overline{N_j}$ & Minimum/maximum number of deliveries offered by contract $j$\\
		$F_{j}$ &    Frequency of deliveries in contract $j$ [hours]               \\
		$O^{+}_{j}, O^{-}_{j} $ &   Maximum up-scaling/down-scaling option offered in contract $j$     [pu]\\
		$\overline{\Delta}$ & Maximum biomass storage level                 [MWt]\\
		$\underline{\Delta}_{t}$ &  Safety storage level of biomass in period   $t$             [MWt]\\
		$\Delta^{\text{F}}$ &    Maximum outflow from biomass storage per period   [MWt/period]\\
		$\Delta^{\text{W}}$ &     Time distance between deliveries to biomass storage [periods]\\
		$\underline{S},\overline{S}$ &  Minimum/maximum thermal storage level                  [MWt]\\
		$S^{\text{F}}$ &     Maximum in/outflow to/from thermal storage per period [MWt/period]\\
		$\underline{P},\overline{P}$ &     Minimum/maximum  production of CHP plant per period [MWe/period]\\
		$\overline{Q^{\text{CHP}}}$ &  Maximum heat production  of CHP plant per period              [MWt/period]\\
		$E ^{\text{CHP}}_{P}$ & Electric efficiency of the CHP plant       [pu]\\
		$E ^{\text{CHP}}_{Q}$ &  Heat efficiency of the CHP plant       [pu]\\
		$E ^{\text{B}}$ &  Calorific value of the biomass      [MWt/tonnes]\\
		$\Theta $ &     Fraction of power reduction         \\
		$\Xi  $ &      Maximum heat to power ratio            \\
		$M^{U}, M^{D} $ &        Minimum up time / down time of CHP plant                       [periods]\\
		$R^{U}, R^{D} $ &        Ramp-up and ramp-down limits of CHP plant                              [MWe/period]\\
		$\overline{Q^{\text{AUX}}}$ &  Maximum heat production of auxiliary boiler per period               [MWt/period]\\
		$E ^{\text{AUX}}$ &  Auxiliary boiler efficiency [pu]\\
		$ P^{\text{B}}$ &  Target percentage of heat produced by biomass [pu] \\
		$\phi^{\text{Sto}}$ & Penalty for excess of storage at the end of time horizon [\euro]\\
		$\phi^{\text{Miss}}$ & Penalty for missed heat demand [\euro] \\
		$\phi^{\text{BM}}$ & Penalty to fail the minimum required heat demand by biomass [\euro]\\
		$\psi_{t}$ & Small incentive for concentrating biomass options in period $t$  [\euro]\\\bottomrule
	\end{tabular}
	\label{tab:parameters}
\end{table}

Biomass from the storage is used by the CHP plant to produce power and heat. The production of both is limited to the feasible production region of an extraction condensing unit depicted with the relevant parameters $\Theta$ and $\Xi$ \cite{zugno2014robust} in Figure \ref{fig:chp_feasibleregion}.  The efficiency of a conversion from biomass to power and heat is denoted by $E^{\text{CHP}}_P$ and $E^{\text{CHP}}_Q$, respectively.
From one hour to the next, the power production of the CHP can be ramped up or down but only in the limits of the parameters $R^U$ and $R^D$. If the unit is started up or shut down it has to be in that state for at least $M^U$ or $M^D$ time periods. Starting up and shutting down is priced with $C^{SU}$ and $C^{SD}$, respectively. The operation of the CHP itself has a cost of $C^{\text{CHP}}$. The power produced is sold on the electricity market and the profit depends on the market price $L^E_{t,\omega}$ in scenario $\omega$. In Denmark, the production of electricity by biomass is supported with an incentive of $I$, while the production of electricity with any fuel is taxed with $T^{\text{EP}}$. Thus, the overall cost $L_{t,\omega}$ is given by $L_{t,\omega} = T^{\text{EP}} - I - L^E_{t,\omega}$, where negative values of $L_{t,\omega}$ are profits. 

The auxiliary boiler has a maximum capacity of $\overline{Q^{\text{AUX}}}$ with an efficiency of $E^{\text{AUX}}$. The operational costs $C^{\text{AUX}}_{t,\omega}$ of the boiler consists of several components and is dependent on the scenario $\omega$ due to the uncertain fuel (e.g. gas or electricity) spot price $C^{\text{F}}_{t,\omega}$. Further components are the operation and maintenance costs $C^{\text{O\&M}}_{AUX}$, taxes  $T^{\text{AUX}}$ and CO$_2$ taxes $T^{\text{CO}_2}$. Thus, the overall operational costs are given by $C^{\text{AUX}}_{t,\omega} = C^{\text{F}}_{t,\omega} + C^{\text{O\&M}}_{AUX} + T^{\text{AUX}} + T^{\text{CO}_2}$.

Both units can feed the thermal storage. In the beginning of the planning horizon (period 0), the heat tank has a given level and the level has to be always between $\underline{S}$ and $\overline{S}$. The in-/outflow per period is limited to $S^F$.

The producer is obliged to fulfill the heat demand in the district heating network $D_{t,\omega}$, which is modeled in scenarios $\omega$.  Due to regulations, the heat production based on biomass is aimed at covering at least $P^B$ percent of the total demand. The probability of scenario $\omega$ is given by $\pi_\omega$. To sum the uncertain parameters up, a scenario $\omega$  resembles the heat demand $D_{t,\omega}$, the electricity price $L^E_{t,\omega}$ and the fuel spot price for the auxiliary boiler $C^{\text{AUX}}_{t,\omega}$.

The overall objective of the solution approach is to select the portfolio of biomass contracts and their configurations that minimizes the cost while fulfilling the heat demand taking the technical characteristics of the plant into account. In this paper, we consider a planning horizon of one year ranging from summer to summer as it is done in practice. Thus, the heating seasons lies in the middle of the planning horizon. However, in general the method can be used with any length of the planning horizon starting and ending at an arbitrary point in time during the year.

\section{Two-phase solution approach}\label{Solutionapproach}

The time scales in the above mentioned planning problem have a broad range. As the contracts are often agreed for up to one year, this results in a medium-term planning problem. However, many technical characteristics of the CHP unit and the electricity market relate to an hourly level.  Additionally, the production does not need to be scheduled more than one week in advance, because then information especially regarding the heat demand gets more accurate. Therefore, we divide the overall planning problem into two-phases:

	 \textbf{Biomass contract selection:} This model decides which suppliers should be contracted for the next year and which amount of biomass they should deliver (including options). The model is based on heat demand scenarios and includes the production by the CHP plant and auxiliary boiler on a weekly time scale excluding ramping and unit commitment decisions. For this long planning horizon the electricity and fuel prices are approximated by an expected value, because the prices are very volatile and hard to predict for a long time horizon. The thermal heat storage is excluded from this model, because it is not reasonable to model the flows on a weekly scale due the small size of those storages. Set $\mathcal{T}$ represents weekly periods in this model. The mathematical formulation is presented in Section \ref{sec:biomass_model}. 
	
	 \textbf{Operational planning problem:} Here the input of biomass is fixed based on the contracts selected in phase 1, but the amounts of contracts with agreed options can still be altered. The model is solved week-by-week taking the input from the previous week into account (storage levels, status of the unit) and decides on the actual production of the CHP plant and auxiliary boiler on an hourly basis incorporating technical requirements and scenario-based price and demand information. Set $\mathcal{T}$ represents hourly periods in this model. The model formulation is described in Section \ref{sec:operational_model}.

Based on the scenario-based representation of the uncertain parameters, both models are two-stage stochastic programming model formulations. 
The division of the planning problem into two phases not only reduces the complexity of the problem, but also resembles the planning process in practice in a more accurate way. Furthermore, solving the operational planning problem week-by-week enables us to make use of more recent information to update the scenarios for the next week. We do not consider an integrated problem for the entire year in an hourly resolution because the addition of such precise information can negatively affect the solution of the problem towards the real realization of the uncertainty due to forecasting inaccuracies. Furthermore, preliminary experiments showed that the large number of integer variables makes the problem computationally hard and not solvable in a reasonable amount of time.

\subsection{Biomass contract selection} \label{sec:biomass_model}

The following model represents the biomass contract selection in phase 1. The model has a weekly time-scale, therefore, the set $\mathcal{T}$ consists of weeks. The relevant parameters like capacities and flow restrictions of the units and storage are scaled up to weekly values accordingly.

\begin{table}[t]
	\caption{Variables}
	\footnotesize
	\begin{tabular}{lp{0.7\columnwidth}}\toprule
		$u_{j} \in \lbrace 0,1 \rbrace$ &   Equals 1, if contract $j$ is used, 0 otherwise      \\
		$d_{j,t} \in \mathbb{N}_{0}$ & Number of deliveries by contract $j$ in period $t$\\
		$\hat{d}_{j,t} \in \lbrace 0,1 \rbrace$ & Equals 1, if contract $j$ delivers in period $t$, 0 otherwise\\
		$b^{}_{j,t} \in \mathbb{R}^+_{0}$ &   Amount of biomass contracted in contract $j$ for period $t$                          [tonnes]\\
		$b^{+}_{j,t}\in \mathbb{R}^+_{0}$ &   Up-scaling option contracted in contract $j$ for period $t$           [tonnes]\\
		$b^{-}_{j,t}\in \mathbb{R}^+_{0}$ &  Down-scaling option contracted in contract $j$ for period $t$            [tonnes]\\
		$\overline{b^{+}}_{j,t,\omega}\in \mathbb{R}^+_{0}$ & Actual amount used of up-scaling option in contract $j$   [tonnes]\\
		$\overline{b^{-}}_{j,t,\omega}\in \mathbb{R}^+_{0}$ & Actual amount used of down-scaling option in contract $j$ [tonnes]\\
		$\delta_{t,\omega}\in \mathbb{R}^+_{0}$ &    Biomass storage level            [MWt]\\
		$\delta^{+}_{t,\omega}\in \mathbb{R}^+_{0}$ &  Inflow to biomass storage     [MWt/period]\\
		$\delta^{{-}}_{t,\omega}\in \mathbb{R}^+_{0}$ &  Outflow from biomass storage  [MWt/period]\\
		$s_{t,\omega}\in \mathbb{R}^+_{0}$ &    Thermal storage level                [MWt]\\
		$s^{{+}}_{t,\omega}\in \mathbb{R}^+_{0}$ & Inflow to thermal storage         [MWt/period]\\
		$s^{{-}}_{t,\omega}\in \mathbb{R}^+_{0}$ & Outflow from thermal storage      [MWt/period]\\
		$x_{t,\omega} \in \lbrace 0,1 \rbrace$ &  Equals 1, if CHP plant is on in period $t$, 0 otherwise    \\       
		$y_{t,\omega} \in \lbrace 0,1 \rbrace$ & Equals 1, if CHP plant is started up in period $t$, 0 otherwise   \\  
		$z_{t,\omega} \in \lbrace 0,1 \rbrace$ & Equals 1, if CHP plant is shut down in period $t$, 0 otherwise   \\  
		$p_{t,\omega}\in \mathbb{R}^+_{0}$ & Power production by CHP  [MWe/period]\\
		$q^{\text{CHP}}_{t,\omega}\in \mathbb{R}^+_{0}$ &     Total heat production by CHP                                [MWt/period]\\
		$q^{\text{CHP,N}}_{t,\omega}\in \mathbb{R}^+_{0}$ &   Heat from CHP flowing to DH    [MWt/period]\\
		$q^{\text{CHP,S}}_{t,\omega}\in \mathbb{R}^+_{0}$ &  Heat from CHP to thermal storage              [MWt/period]\\
		
		$q^{\text{AUX}}_{t,\omega} \in \mathbb{R}^+_{0} $ &  Total heat production by auxiliary boiler                                [MWt/period]\\
		$q^{\text{AUX,N}}_{t,\omega}\in \mathbb{R}^+_{0}$ &    Heat from auxiliary boiler to DH  [MWt/period]\\
		$q^{\text{AUX,S}}_{t,\omega}\in \mathbb{R}^+_{0}$ & Heat from auxiliary boiler to thermal storage              [MWt/period]\\
		
		$q^{\text{Miss}}_{t,\omega}\in \mathbb{R}^+_{0}$ & Missed heat demand                  [MWt/period]\\
		$q_{t,\omega}^{\text{BM}}\in \mathbb{R}^+_{0}$ & Required amount of heat not supplied with biomass [MWt/period]\\
		$\delta_{t,\omega}^{\text{EX}}\in \mathbb{R}^+_{0}$ & Amount of biomass above storage capacity [MWt]\\
		$\delta_{\omega}^{\text{T}}\in \mathbb{R}^+_{0}$ & Amount of biomass in excess at the end of the time horizon [MWt]\\\bottomrule
	\end{tabular}
	\label{tab:variables}
\end{table}

The first-stage decision variables in this model decide on the contracts to be selected ($u_j$) as well as the number of deliveries in each week ($d_{j,t}$) and amounts ($b_{j,t}$) including up- ($b^{+}_{j,t}$) and down-scaling ($b^{-}_{j,t}$) options for each contract $j \in \mathcal{J}$ and period $t \in \mathcal{T}$. Based on the second-stage variables, these amounts can be altered with the variables $\overline{b^{+}}_{j,t,\omega}$ and $\overline{b^{-}}_{j,t,\omega}$ within the limits of the selected options in the first-stage. Further second-stage variables relate to the biomass storage level ($\delta_{t,\omega}$) as well as heat ($q^{\text{CHP}}_{t,\omega},q^{\text{AUX}}_{t,\omega}$) and power production ($p_{t,\omega}$). An overview of the variables and their domains is given in Table \ref{tab:variables}.

{\allowdisplaybreaks
\begin{subequations} \label{eq1:obj}
	\begin{align}
	&{\text{min}} \negthickspace \sum_{t\ \in \mathcal{T}} \negthickspace\Bigg[ \sum_{j \in \mathcal{J}} \bigg(C^{\text{B}}_{j}b^{}_{j,t}+C^{\text{B}+}_{j}b^{+}_{j,t}+C^{\text{B}-}_{j}b^{-}_{j,t}
	+ \negthickspace\negthickspace\sum_{\omega \in \Omega}\negthickspace \pi_{\omega} C^{\text{B}}_{j}(\overline{b^{+}}_{j,t,\omega}-\overline{b^{-}}_{j,t,\omega})\bigg)\label{eq1:obj_contractcosts}\\
	& +\negthickspace\sum_{\omega \in \Omega} \pi_{\omega} \bigg( C^{\text{CHP}}\Big[p_{t,\omega}-\Theta q^{\text{CHP}}_{t,\omega}\Big]+\widehat{L_{t}}p_{t,\omega}+  \widehat{C^{\text{AUX}}_{t}}\frac{q^{\text{AUX}}_{t,\omega}}{E^{\text{AUX}}}+ C^\text{I} \delta_{t,\omega}\bigg)\Bigg]\label{eq1:obj_productioncosts}\\
	&+\negthickspace\sum_{\omega \in \Omega}\pi_{\omega} \bigg( \phi^{\text{Sto}}\delta^{\text{T}}_{\omega}  + \sum_{t \in \mathcal{T}} \big( \phi^{\text{Miss}}q^{\text{Miss}}_{t,\omega} + \phi^{\text{BM}} q^{\text{BM}}_{t,\omega}\big) \bigg) \negthickspace- \negthickspace \sum_{t\ \in \mathcal{T}} \sum_{j \in \mathcal{J}} \psi_t (b^{+}_{j,t}+b^{-}_{j,t})
	\label{eq1:obj_penalty}
	\end{align}
\end{subequations}}
The objective function \eqref{eq1:obj} minimizes the expected cost of the biomass contract selection. The first part \eqref{eq1:obj_contractcosts} contains the costs related to the biomass supply and the contract selection. In \eqref{eq1:obj_productioncosts}, operational costs of the system, profits from electricity sales and inventory costs for biomass are modeled. Note, that the parameters 
$\widehat{L_{t}}$ and $\widehat{C^{\text{AUX}}_{t}}$ are expected values for this tactical problem. The third part \eqref{eq1:obj_penalty} represents penalty and virtual costs. First, we penalize leftover biomass at the end of the planning period ($\phi^{\text{Sto}}\delta^{\text{T}}_{\omega}$), since we try to empty the storage at the end of the year. Second,   missed heat-demand ($\phi^{\text{Miss}}q^{\text{Miss}}_{t,\omega}$) is penalized. Finally, we add a penalty for failing to meet the minimum share of heat production by biomass ($\phi^{\text{BM}} q^{\text{BM}}_{t,\omega}$). The second sum incentivizes the use options in certain periods with a very small profit ($\psi_t$). This allows to concentrate options in periods with a high variance in scenarios. In preliminary experiments it turned out, that there are equally good solutions as the price for options is the same over the year. When the inventory costs are low, the options and amounts can be shifted without deteriorating the objective. Therefore, we introduce this incentive to prioritize weeks with a high variance in  demand.
{\allowdisplaybreaks
	\begin{align}
	\underline{N_{j}}u_{j} &\leq \sum_{t \in \mathcal{T}} d_{j,t} \leq \overline{N_{j}} u_{j} & \forall j \in \mathcal{J} \label{eq1:activation_contract}\\ 
	\sum_{\tau = t- \min{\lbrace \lfloor \frac{F_{j}}{168}} \rfloor, 1 \rbrace}^{t}d_{j,\tau} &\leq \max{\Bigg\lbrace \frac{168}{F_j}, 1 \Bigg\rbrace} & \forall j \in \mathcal{J},\forall t \in \mathcal{T} \label{eq1:frequency}\\
	b^{}_{j,t}+b^{+}_{j,t} &\leq \overline{B_{j}}d_{j,t} & \forall j \in \mathcal{J}, \forall t \in \mathcal{T} \label{eq1:ub_contract}\\
	b^{}_{j,t}-b^{-}_{j,t} &\geq \underline{B_{j}}d_{j,t} & \forall j \in \mathcal{J}, \forall t \in \mathcal{T} \label{eq1:lb_contract} \\
	b^{+}_{j,t} &\leq O^{+}_{j}b^{}_{j,t} & \forall j \in \mathcal{J}, \forall t \in \mathcal{T}  \label{eq1:ub_option}\\
	b^{-}_{j,t} &\leq O^{-}_{j}b^{}_{j,t}  & \forall j \in \mathcal{J}, \forall t \in \mathcal{T} \label{eq1:lb_option}\\
	\overline{b^{+}}_{j,t,\omega} &\leq b^{+}_{j,t} & \hspace{-1.5cm} \forall j \in \mathcal{J},\forall t \in \mathcal{T}, \forall \omega \in \Omega \label{eq1:ub_option_actual}\\
	\overline{b^{-}}_{j,t,\omega} &\leq b^{-}_{j,t} & \hspace{-1.5cm} \forall j \in \mathcal{J},\forall t \in \mathcal{T}, \forall \omega \in \Omega  \label{eq1:lb_option_actual}
	\end{align}
	Constraints \eqref{eq1:activation_contract} to \eqref{eq1:lb_option_actual} model the selection of biomass contracts. In constraints \eqref{eq1:activation_contract} the number of deliveries is restricted by the contract limits. Constraint \eqref{eq1:frequency} restricts the number of deliveries per week to a maximum according to the frequency of the contract. The left-hand side sums over several weeks, if the minimum time between visits $F_j$ is longer than one week (168 hours). The right-hand side determines the maximum number of deliveries in that period with at least one delivery or more if the time difference is less than 168 hours.  The total amount including up- and down-scaling options is limited by constraints \eqref{eq1:ub_contract} and \eqref{eq1:lb_contract} and the use of options in constraints \eqref{eq1:ub_option} and \eqref{eq1:lb_option}. In constraints \eqref{eq1:ub_option_actual} and \eqref{eq1:lb_option_actual}, it is ensured that the second-stage alterations respect the first-stage decisions.
	\begin{align}
	\underline{\Delta}_t  &\leq \delta_{t,\omega} \leq \overline{\Delta} & \forall t \in \mathcal{T}, \forall \omega \in \Omega \label{eq1:storage_limits} \\
	\delta_{t,\omega}  &= \delta_{t-1,\omega}+\delta^{{+}}_{t,\omega}-\delta^{{-}}_{t,\omega} & \forall t \in \mathcal{T}, \forall \omega \in \Omega \label{eq1:storage_level}\\
	\delta^{{+}}_{t,\omega}&= \sum_{j \in \mathcal{J}} \Big(b_{j,t}+ \overline{b^{+}}_{j,t,\omega}- \overline{b^{-}}_{j,t,\omega} \Big) \cdot E^{\text{B}}  &\forall t \in \mathcal{T}, \forall \omega \in \Omega \label{eq1:storage_inflow}\\ 
	\delta^{{+}}_{t,\omega}&\leq \Delta^{\text{F}} & \forall t \in \mathcal{T}, \forall \omega \in \Omega \label{eq1:storage_flow}\\
	\delta_{|\mathcal{T}|,\omega} &\leq \delta_{0,\omega}+\delta^{\text{T}}_{\omega} &  \forall \omega \in \Omega \label{eq1:storage_end} 
	\end{align}	
	The biomass storage is modeled by constraints \eqref{eq1:storage_limits} to \eqref{eq1:storage_end}. The model ensures that the storage level is kept within the limits \eqref{eq1:storage_limits} and calculated correctly based on the previous level and in- and outflows \eqref{eq1:storage_level}. The initial storage level is given by $\delta_{0,\omega}$, which is the same for all scenarios. The inflow from supplier deliveries is calculated in constraints \eqref{eq1:storage_inflow}, where the incoming biomass is converted from tonnes to MWht using the calorific value of the biomass $E^{\text{B}}$ and the outflow is restricted by constraints \eqref{eq1:storage_flow}. Finally, the storage level at the end of the planning horizon is determined in \eqref{eq1:storage_end}  for penalty cost calculations.
	\begin{align}	
		\delta^{{-}}_{t,\omega} &=\frac{ p_{t,\omega}}{E ^{\text{CHP}}_{P}} - \Theta \cdot \frac{q^{\text{CHP}}_{t,\omega}}{E ^{\text{CHP}}_{Q} }& \forall t \in \mathcal{T}, \forall \omega \in \Omega \label{eq1:biomass_consumption}\\
	\underline{P_{}} &\leq p_{t,\omega}-\Theta \cdot q^{\text{CHP}}_{t,\omega}  \leq \overline{P_{}} &\forall t \in \mathcal{T}, \forall \omega \in \Omega \label{eq1:chp_feasreg1} \\
	\Xi  \cdot q^{\text{CHP}}_{t,\omega}  & \leq  p_{t,\omega} & \forall t \in \mathcal{T}, \forall \omega \in \Omega \label{eq1:chp_feasreg2}\\
	q^{\text{CHP}}_{t,\omega}  & \leq \overline{Q^{\text{CHP}}}  & \forall t \in \mathcal{T}, \forall \omega \in \Omega \label{eq1:chp_feasreg3}\\
	q^{\text{AUX}}_{t,\omega} & \leq \overline{Q^{\text{AUX}}_{}}  &\forall t \in \mathcal{T}, \forall \omega \in \Omega  \label{eq1:gb_feasreg}
	\end{align}	
	The production capacities of the CHP plant and auxiliary boiler are enforced by constraints \eqref{eq1:biomass_consumption} to \eqref{eq1:gb_feasreg}. In \eqref{eq1:biomass_consumption} the consumption of biomass from the storage for CHP production is determined based on the corresponding efficiency. The feasible region of the CHP, which was previously presented in Figure \ref{fig:chp_feasibleregion}, is modeled by constraints \eqref{eq1:chp_feasreg1} to \eqref{eq1:chp_feasreg3} and  limits of the auxiliary boiler in \eqref{eq1:gb_feasreg}.
	\begin{align}	
	D_{t,\omega} &= q^{\text{CHP}}_{t,\omega} + q^{\text{AUX}}_{t,\omega} + q^{\text{Miss}}_{t,\omega} &  \forall t \in \mathcal{T}, \omega \in \Omega \label{eq1:heatdemand}\\
	q^{\text{CHP}}_{t,\omega} & \ge P^{\text{B}}\cdot D_{t,\omega} - q^{\text{BM}}_{t,\omega} &  \forall t \in \mathcal{T}, \omega \in \Omega \label{eq1:minbiomassheat}
	\end{align}
	Finally, the heat demand is ensured in constraint \eqref{eq1:heatdemand} while at least $P^{\text{B}}$ percent  per week have to be produced by biomass otherwise causing penalty costs \eqref{eq1:minbiomassheat}. }


\subsection{Operational planning} \label{sec:operational_model}
The operational planning model relates to the second phase of the solution approach. For the overall solution approach, the model is solved consecutively week-by-week with a receding horizon to determine the production schedule and to adjust the biomass deliveries, if possible.  Therefore, the planning horizon is $|\mathcal{W}|$ weeks with an hourly resolution. The week in focus is $\mathcal{W}_{1}$ and the remaining weeks $\mathcal{W}_{2}$ to $\mathcal{W}_{\mathrm{w}}$ are used in the receding horizon to already include predictions for future periods. Thus, the decisions for weeks $\mathcal{W}_{2}$ to $\mathcal{W}_{\mathrm{w}}$ can be altered again later, when the respective week comes in focus. Set $\mathcal{T}$ consists of all hours in the planning horizon, whereas  $\mathcal{T}_\mathrm{w}$ relates to the hours in specific week $\mathrm{w} \in \mathcal{W}$.

The decision variables for this model decide the amount ($b_{j,t,\omega}$) and up- and down- scaling biomass  ($\overline{b^{+}}_{j,t,\omega}$ and $\overline{b^{-}}_{j,t\omega}$) and the actual delivery times ($\hat{d}_{j,t,\omega}$) for the deliveries of contract $j$. Further variables are related to the biomass storage level ($\delta_{t,\omega}$), the thermal storage level ($s_{t,\omega}$), the heat and power production ($q^{\text{CHP}}_{t,\omega}$, $q^{\text{AUX}}_{t,\omega}$ and $p_{t,\omega}$) and the commitment status of the CHP plant ($x_{t,\omega}$, $y_{t,\omega}$ and $z_{t,\omega}$). The variables are also included in Table \ref{tab:variables}.

Because the first week of the receding horizon is the week in focus, the first-stage decisions of the stochastic program are the delivery times and amounts $\hat{d}_{j,t,\omega}$, $b_{j,t,\omega}$, $\overline{b^{+}}_{j,t,\omega}$ and $\overline{b^{-}}_{j,t\omega}$  for periods $t$ in the first week $\mathcal{T}_1$. For all other weeks, the decisions can be revised later and therefore are second-stage decisions. To ensure non-anticipativity, we include specific constraints.

The selection of biomass contracts and amounts are input parameters to this model (given in Table \ref{tab:input}) and determined by the biomass contract selection model in phase 1. Set $\mathcal{J}$ is reduced to only selected contracts for the corresponding week to limit the number of variables.

Furthermore, the storage levels and unit status of the preceding week are set as initial values. For example, the initial biomass storage level $\delta_{0,\omega}$ in the current week  equals the storage level in the last period of previous week.

\begin{subequations}\label{eq2:obj}
	\begin{align}
	&{\text{min}} \sum_{\mathrm{w} \in \mathcal{W}}  \sum_{j \in \mathcal{J}} \bigg(C^{\text{B}}_{j}\mathbf{B_{j,\mathrm{w}}}+C^{\text{B}+}_{j}\mathbf{B^+_{j,\mathrm{w}}}+C^{\text{B}-}_{j}\mathbf{B^-_{j,\mathrm{w}}}  \sum_{t \in \mathcal{T}_\mathrm{w}}
	C^{\text{B}}_{j}(\overline{b^{+}}_{j,t}-\overline{b^{-}}_{j,t})\bigg)\label{eq2:obj_biomass}\\
	& + \sum_{t \in \mathcal{T}}    \sum_{\omega \in \Omega}  \pi_{\omega} \bigg( C^{\text{CHP}}\big(p_{t,\omega}-\Theta q^{\text{CHP}}_{t,\omega}\big) \negthickspace- \negthickspace{{L_{t,\omega}}}p_{t,\omega}
	+ C^{\text{SU}}y_{t,\omega}  + C^{\text{SD}}z_{t,\omega}\bigg) \label{eq2:obj_operational} \\
	& + \sum_{t \in \mathcal{T}}    \sum_{\omega \in \Omega} \pi_{\omega} \big({C^{\text{AUX}}_{t,\omega}}\frac{q^{\text{AUX}}_{t,\omega}}{E^{\text{AUX}}} + C^\text{I}\delta_{t,\omega} \big) \negthickspace+ \negthickspace\sum_{t \in \mathcal{T}}    \sum_{\omega \in \Omega} \pi_{\omega} \big(\phi^{\text{Sto}}\delta^{\text{EX}}_{t,\omega} +  \phi^{\text{Miss}}q^{\text{Miss}}_{t,\omega} \big)  \label{eq2:obj_penalty}
	\end{align}
\end{subequations}

\begin{table}
	\caption{Input parameters from biomass contract selection}
	\footnotesize
	\begin{tabular}{ll}\toprule
		$\mathbf{U_{j,\mathrm{w}}} \in \mathbb{N}_{0}$& Number of deliveries of contract $j$ in week $\mathrm{w}$ \\
		$\mathbf{B_{j,\mathrm{w}}} \in \mathbb{R}^+_{0}$& Contracted delivery amount of contract $j$ in week $\mathrm{w}$ \\
		$\mathbf{B^+_{j,\mathrm{w}}} \in \mathbb{R}^+_{0}$& Contracted up-scaling of delivery amount of contract $j$ in week $\mathrm{w}$ \\
		$\mathbf{B^-_{j,\mathrm{w}}} \in \mathbb{R}^+_{0}$& Contracted down-scaling of delivery amount of contract $j$ in  week $\mathrm{w}$ \\\bottomrule
	\end{tabular} 
	\label{tab:input}
\end{table}

As in the biomass contract selection model, the objective function \eqref{eq2:obj} minimizes the expected costs composed of biomass contract costs \eqref{eq2:obj_biomass}, operational for the CHP \eqref{eq2:obj_operational}, operational costs for the auxiliary and the biomass storage \eqref{eq2:obj_penalty}, and penalty  costs \eqref{eq2:obj_penalty}. However, the following changes have to be made. First, the profit for electricity sales (${{L_{t,\omega}}}$) and operational costs for the auxiliary boiler (${C^{\text{AUX}}_{t,\omega}}$) depend on scenarios \eqref{eq2:obj_operational}. Second, the operational cost \eqref{eq2:obj_operational} now includes costs for starting up and shutting down the CHP plant. Third, the term \eqref{eq2:obj_penalty} penalizes unfulfilled heat demands and exceeding the biomass storage capacity. Note that to resemble the total weekly cost of the system, we keep the constant term $C^{\text{B}}_{j}\mathbf{B_{j,\mathrm{w}}}+C^{\text{B}+}_{j}\mathbf{B^+_{j,\mathrm{w}}}+C^{\text{B}-}_{j}\mathbf{B^-_{j,\mathrm{w}}}$ in \eqref{eq2:obj_biomass}. 
{\allowdisplaybreaks
	\begin{align}
	\sum_{t \in \mathcal{T_\mathrm{w}}} \hat{d}_{j,t,\omega} &= \mathbf{U_{j,\mathrm{w}}} && \forall j \in \mathcal{J}, \forall \mathrm{w} \in \mathcal{W}, \forall \omega \in \Omega\label{eq2:delivery}\\
	\sum_{t \in \mathcal{T_\mathrm{w}}} b_{j,t,\omega} &= \mathbf{B_{j,\mathrm{w}}} && \forall j \in \mathcal{J}, \forall \mathrm{w} \in \mathcal{W}, \forall \omega \in \Omega\label{eq2:deliveryamount}\\
	\sum_{t \in \mathcal{T_\mathrm{w}}}\overline{b^{+}}_{j,t,\omega} &\leq \mathbf{B^{+}_{j,\mathrm{w}}}  && \forall j \in \mathcal{J},\forall \mathrm{w} \in \mathcal{W}, \forall \omega \in \Omega\label{eq2:upscaling} \\
	\sum_{t \in \mathcal{T_\mathrm{w}}}\overline{b^{-}}_{j,t,\omega} &\leq \mathbf{B^{-}_{j,\mathrm{w}}} && \forall j \in \mathcal{J},\forall \mathrm{w} \in \mathcal{W}, \forall \omega \in \Omega\label{eq2:downscaling}  \\
	b_{j,t,\omega} + \overline{b^+}_{j,t,\omega} &\le \overline{B}_j \hat{d}_{j,t,\omega} && \forall j \in \mathcal{J}, \forall t \in \mathcal{T}, \forall \omega \in \Omega\label{eq2:delamountub}\\
	b_{j,t,\omega} - \overline{b^-}_{j,t,\omega} &\ge \underline{B}_j \hat{d}_{j,t,\omega} && \forall j \in \mathcal{J}, \forall t \in \mathcal{T}, \forall \omega \in \Omega\label{eq2:delamountlb}\\
	\sum_{\tau=t-F_j}^{t} \hat{d}_{j,\tau,\omega} & \le 1 && \forall j \in \mathcal{J}, \forall t \in \mathcal{T}, \forall \omega \in \Omega\label{eq2:delfreq}\\
	\sum_{j \in \mathcal{J}}\sum_{\tau=t}^{t+{\Delta^\text{W}}}\hat{d}_{j,\tau,\omega} &\leq 1 && \forall t \in \mathcal{T},\forall \omega \in \Omega  \label{eq2:wait}
	\end{align}	
	The biomass deliveries are handled in constraints \eqref{eq2:delivery} to \eqref{eq2:wait}. If deliveries were scheduled for the weeks in the planning horizon by phase 1, the operation model decides on the actual delivery times during the week \eqref{eq2:delivery}. The weekly contracted amount is split on the deliveries in constraints \eqref{eq2:deliveryamount}. The delivery amount can be altered in the given limits of the options (constraints \eqref{eq2:upscaling} and \eqref{eq2:downscaling}), but the total amount must be within the limits of the contract (constraints \eqref{eq2:delamountub} and \eqref{eq2:delamountlb}). Constraints \eqref{eq2:delfreq} imposes a maximum frequency on the  deliveries associated with each contract, while constraints \eqref{eq2:wait} ensures an elapsed time of at least $\Delta^{\text{W}}$ periods between two deliveries irrespective of the supplier.
	\begin{align}
	\hat{d}_{j,t,\omega} = \hat{d}_{j,t,\omega'},& \quad{b}_{j,t,\omega} = {b}_{j,t,\omega'} && \negthickspace \negthickspace \negthickspace\forall j \in \mathcal{J}, \forall t \in \mathcal{T}_1, \forall \omega,\omega'  \in \Omega, \omega \neq\omega' \label{eq2:nonanticipativityd}\\
	\overline{b^+}_{j,t,\omega} = \overline{b^+}_{j,t,\omega'},& \quad \overline{b^-}_{j,t,\omega} = \overline{b^-}_{j,t,\omega'}  && \negthickspace \negthickspace \negthickspace \forall j \in \mathcal{J}, \forall t \in \mathcal{T}_1, \forall \omega,\omega'  \in \Omega, \omega \neq\omega'\label{eq2:nonanticipativityb+}
	\end{align}
	As the decisions for the biomass delivery in the first week are first-stage decisions of the stochastic program, we have to ensure that they have the same values for each scenario. This is forced by the non-anticipativity constraints 
	\eqref{eq2:nonanticipativityd} to \eqref{eq2:nonanticipativityb+}.
	\begin{align}
	\delta^{{+}}_{t,\omega}& = \sum_{j \in \mathcal{J}} \Big(b_{j,t,\omega}+ \overline{b^{+}}_{j,t,\omega}- \overline{b^{-}}_{j,t,\omega} \Big) \cdot E^{\text{B}}  &\forall t \in \mathcal{T}, \forall \omega \in \Omega \label{eq2:biomass_inflow}\\ 
	\delta_{t,\omega}  &= \delta_{t-1,\omega}+\delta^{{+}}_{t,\omega}-\delta^{{-}}_{t,\omega} & \forall t \in \mathcal{T}, \forall \omega \in \Omega \label{eq2:storage_level}\\
	\delta^{{-}}_{t,\omega}&\leq \Delta^{\text{F}}  & \forall t \in \mathcal{T}, \forall \omega \in \Omega \label{eq2:biomass_maxflow} \\
	\delta_{t,\omega} & \leq \overline{\Delta} + \delta^{\text{EX}}_{t,\omega}& \forall t \in \mathcal{T}, \forall \omega \in \Omega \label{eq2:biomass_capacity_max}\\
	\underline{\Delta}_t  &\leq \delta_{t,\omega} & \hspace{-35pt}   \forall \mathrm{w} \in \lbrace 2,\ldots,|\mathcal{W}|\rbrace, \forall t \in \mathcal{T}_{\mathrm{w}}, \forall \omega \in \Omega \label{eq2:biomass_capacity_safety}\\
	0  &\leq \delta_{t,\omega} &  \forall t \in \mathcal{T}_{1}, \forall \omega \in \Omega \label{eq2:biomass_capacity_week1}
	\end{align}
	The inflow to the biomass storage in each period \eqref{eq2:biomass_inflow} is dependent on the scheduled delivery and adjustments based on the options. The storage level is given by equation \eqref{eq2:storage_level}. The outflow and capacity of the storage is limited in constraints \eqref{eq2:biomass_maxflow} and \eqref{eq2:biomass_capacity_max}, respectively.  The safety storage for biomass is incorporated in constraints \eqref{eq2:biomass_capacity_safety}, but only for future weeks in the receding horizon. In the current week, the storage can be used for production \eqref{eq2:biomass_capacity_week1}.
	\begin{align}
	\delta^{{-}}_{t,\omega} &=\frac{ p_{t,\omega}}{E ^{\text{CHP}}_{P}} - \Theta \cdot \frac{q^{\text{CHP}}_{t,\omega}}{E ^{\text{CHP}}_{Q} }& \forall t \in \mathcal{T}, \forall \omega \in \Omega \label{eq2:biomassconsumption} \\
	\underline{P_{}}\cdot x_{t,\omega}& \leq p_{t,\omega}-\Theta \cdot q^{\text{CHP}}_{t,\omega}  \leq \overline{P_{}}\cdot x_{t,\omega}  & \forall t \in \mathcal{T}, \forall \omega \in \Omega  \label{eq2:chp_feasreg1}\\
	\Xi  \cdot q^{\text{CHP}}_{t,\omega}  & \leq  p_{t,\omega} & \forall t \in \mathcal{T}, \forall \omega \in \Omega \label{eq2:chp_feasreg2}\\
	q^{\text{CHP}}_{t,\omega}  & \leq \overline{Q^{\text{CHP}}} \cdot x_{t,\omega} & \forall t \in \mathcal{T}, \forall \omega \in \Omega \label{eq2:chp_feasreg3}\\
	y_{t,\omega}-z_{t,\omega} & = x_{t,\omega}-x_{t-1,\omega} & \forall t \in \mathcal{T}, \forall \omega \in \Omega \label{eq2:chp_status}\\
	y_{t,\omega}+z_{t,\omega} &\leq 1 & \forall t \in \mathcal{T}, \forall \omega \in \Omega\label{eq2:notstartandstop}\\
	\sum^{t}_{\tau=t-M^U_{}+1} y_{\tau,\omega}& \leq x_{t,\omega}  & \forall t \in \mathcal{T}, \forall \omega \in \Omega \label{eq2:minuptime}\\
	\sum^{t}_{\tau=t-M^D_{}+1} z_{\tau,\omega} &\leq 1-x_{t,\omega} & \forall t \in \mathcal{T}, \forall \omega \in \Omega\label{eq2:mindowntime}\\
	p_{t,\omega}-p_{t-1,\omega} & \leq R^U \cdot x_{t-1,\omega}+  \underline{P_{}} \cdot y_{t-1,\omega} & \forall t \in \mathcal{T}, \forall \omega \in \Omega\label{eq2:rampup}\\
	p_{t,\omega}-p_{t-1,\omega} & \geq -R^D \cdot x_{t,\omega}-  \underline{P_{}} \cdot z_{t,\omega} & \forall t \in \mathcal{T}, \forall \omega \in \Omega\label{eq2:rampdown}
	\end{align}
	Constraints \eqref{eq2:biomassconsumption} to \eqref{eq2:chp_feasreg3} regarding biomass consumption and feasible production region of the CHP unit constraints are similar to constraints \eqref{eq1:biomass_consumption} to \eqref{eq1:chp_feasreg3} for the biomass selection problem. However, here the production depends also on the status of the unit ($x_{t,\omega} = 1$ means the unit is on). The status of the unit is determined by constraints \eqref{eq2:chp_status} to \eqref{eq2:notstartandstop} while constraints \eqref{eq2:minuptime} and \eqref{eq2:mindowntime} ensure minimum up- and down times, respectively. The change of production volume is restricted to the ramping requirements in constraints \eqref{eq2:rampup} and \eqref{eq2:rampdown}. The initial status of the CHP plant depends on the previous week and is given by $x_{0,\omega}$ and $p_{0,\omega}$ as input parameters.
	\begin{align}
	q^{\text{AUX}}_{t,\omega} & \leq \overline{Q^{\text{AUX}}_{}} & \forall t \in \mathcal{T}, \forall \omega \in \Omega \label{eq2:maxgb}\\
	s^{{+}}_{t,\omega} & = q^{\text{CHP,S}}_{t,\omega}+q^{\text{GB,S}}_{t,\omega} & \forall t \in \mathcal{T}, \forall \omega \in \Omega\label{eq2:thermal_inflow}\\
	s_{t,\omega} &=s_{t-1,\omega}+s^{{+}}_{t,\omega}-s^{{-}}_{t,\omega}  & \forall t \in \mathcal{T}, \forall \omega \in \Omega\label{eq2:thermal_level}\\
	\underline{S} & \leq s_{t,\omega} \leq \overline{S}  & \forall t \in \mathcal{T}, \forall \omega \in \Omega\label{eq2:thermal_limits}\\
	s^{{-}}_{t,\omega} & \leq S^{\text{F}}  & \forall t \in \mathcal{T}, \forall \omega \in \Omega\label{eq2:thermal_maxoutflow}\\
	s^{{+}}_{t,\omega} & \leq S^{\text{F}}  & \forall t \in \mathcal{T}, \forall \omega \in \Omega\label{eq2:thermal_maxinflow}\\
	s^{{-}}_{t,\omega} & \leq s_{t-1,\omega}  & \forall t \in \mathcal{T}, \forall \omega \in
	\Omega\label{eq2:thermalshift}\\
	s_{|\mathcal{T_{\mathrm{w}}}|,\omega} & =s_{0,\omega}  & \forall \omega \in \Omega\label{eq2:thermal_end}
	\end{align}
Constraints \eqref{eq2:maxgb} sets the heat production capacity of the auxiliary boiler.
	The heat storage is modeled by constraints \eqref{eq2:thermal_inflow} to \eqref{eq2:thermal_end}. The inflow is determined by the heat from the CHP unit and auxiliary boiler inserted into the storage \eqref{eq2:thermal_inflow}. The current storage level depends on the inflow, outflow and previous level \eqref{eq2:thermal_level} ($s_{0,\omega}$ for the initial value) and has to satisfy the capacity restrictions \eqref{eq2:thermal_limits}. Outflow \eqref{eq2:thermal_maxoutflow} and inflow \eqref{eq2:thermal_maxinflow} are limited and the inflow cannot directly flow out again \eqref{eq2:thermalshift}. To avoid emptying the storage at the end, the initial level $s_{0,\omega}$ must be reached again at the end of the receding horizon \eqref{eq2:thermal_end}. 
	\begin{align}
	q^{\text{CHP}}_{t,\omega} &= q^{\text{CHP,N}}_{t,\omega}+q^{\text{CHP,S}}_{t,\omega} & \forall t \in \mathcal{T}, \forall \omega \in \Omega \label{eq2:heat_chp} \\ 
	q^{\text{AUX}}_{t,\omega} &= q^{\text{AUX,N}}_{t,\omega}+q^{\text{AUX,S}}_{t,\omega}  & \forall t \in \mathcal{T}, \forall \omega \in \Omega \label{eq2:heat_gb}\\
	D_{t,\omega} &=  q^{\text{CHP,N}}_{t,\omega}+ q^{\text{AUX,N}}_{t,\omega}+s^{{-}}_{t,\omega}+q^{\text{Miss}}_{t,\omega} & \forall t \in \mathcal{T}, \forall \omega \in \Omega \label{eq2:heat_demand}
	\end{align}
	The heat production by both units is used for filling the heat storage and covering the demand. Therefore, the production is split up into those two components in constraints \eqref{eq2:heat_chp} and \eqref{eq2:heat_gb}.  For fulfilling the heat demand, heat directly fed to the district heating network and heat from the thermal storage  is used \eqref{eq2:heat_demand}. Any shortfall of heat is penalized in the objective function. 
}

\subsection{Overall solution approach}
For the overall solution approach, the above mentioned stochastic programming models are combined. To solve the planning problem for one year, we need to perform the steps shown in Algorithm \ref{alg:approach}. First, the contract selection takes place (line 1). Afterward, this decision is transfered to the weekly planning (line 3). The scenario generation and solution of the operational problem is carried once every week for the next week (line 2 to 6).

\begin{algorithm}[t]
	\footnotesize
	\caption{Two-phase solution approach}
	\begin{algorithmic}[1]
		\State	Solve the biomass contract selection model \eqref{eq1:obj}-\eqref{eq1:minbiomassheat}
		\For{each week in the overall planning horizon}
		\State Select the corresponding contract decisions from line 1 and set limits 
		\State Generate scenarios for the current receding horizon
		\State Solve the operational planning model \eqref{eq2:obj}-\eqref{eq2:heat_demand}
		\EndFor
	\end{algorithmic}
	\label{alg:approach}
\end{algorithm}


\section{Case studies}\label{casestudies} 
In the following we analyze two case studies for different municipalities in Denmark, named A and B, that are connected to the Aarhus district heating network.  The planning horizon we consider in the numerical results in Section \ref{sec:results} is 1st of June 2016 to  31st of May 2017.

\subsection{Technical data}
The heat demand data in the district heating networks is obtained from \cite{Dataadga99:online}, NordPools' hourly electricity prices for DK1 zone from \cite{Dataabou20:online} and daily natural gas prices from \cite{Powernex20:online}. Extreme outlier values in electricity prices are limited to a maximum or minimum of four standard deviations from the mean.

The technical parameters for the CHP and auxiliary units as well as the operation costs are based on \cite{nielsen2016economic,soysal2016electric,levitt2014cost} and \cite{197363:online} and shown in Tables \ref{tab:technicalchp}, \ref{tab:technicalaux} and \ref{tab:cost}. Both systems comprise a CHP unit and one auxiliary boiler. Municipality A uses a gas boiler in addition to the CHP, while municipality B uses an electric boiler. The biomass storage minimum level $\underline{\Delta}_t$ is divided in two values. In weeks 20 - 45 (i.e. in the heating season), we have a higher minimum level as in the remaining weeks of the year. The penalty costs for both case are the same and set to $\Phi^{\text{Sto}} =	1000       , \Phi^{\text{Miss}}=10000$ and  $\Phi^{\text{BM}} = 5000$.

The parameters of the biomass contracts data are given in Table \ref{tab:contracts}, where they are organized from \textit{fixed} contracts at the top of the table and gradually going down to more \textit{flexible} contracts.  Both cases use the same set of contracts. 

The very small incentive $\psi_{t}$ for using options preferably in periods with a high variance in heat demand scenarios is calculated as follows. As this is a weekly value for the biomass contract selection phase only, we consider the weekly heat demand scenarios in phase 1. This data is known before solving the model and therefore we can use the scenario information to calculate the incentive. We order the weeks $t$ in descending order of difference in heat demands in the scenarios, i.e, $max_{\omega \in \Omega} \lbrace D_{t,\omega} \rbrace - min_{\omega \in \Omega} \lbrace D_{t,\omega} \rbrace$. The week with the largest difference gets the highest incentive of $5.2$. We reduce the incentive every week by $0.1$ resulting in an incentive of $0.1$ for the week with the smallest difference. These values are far less than the cost of the biomass options themselves and, therefore, have barely influence on the amounts contracted in options but only on the weeks where they are placed. Note that this incentive is not part of the evaluation in Section \ref{sec:results} as it is only in the biomass contract selection and the cost  are based on the operational planning.

\begin{table*}
	\caption{Technical parameters of the CHP unit}

	\footnotesize
	\centering
	\begin{tabular}{llllllllllll}\toprule		
		&    $\overline{P}$               &    $\underline{P}$              &    $\overline{Q^{\text{CHP}}}$  &    $\Theta$                     &    $\Xi$                        &    $R^{U}$                      &    $R^{D}$                      &    $E^{\text{CHP}}_{P}$         &      $E^{\text{CHP}}_{Q}$       &    $M^{U}$                      &    $M^{D}$                      \\\midrule
		
		A &      13.24 &        3.8 &       20.8 &      -0.18 &       0.55 &        3.7 &        3.7 &       0.62 &       0.31 &          6 &          4 \\
		
		B &      35.18 &       5.72 &      47.28 &      -0.12 &       0.64 &        4.6 &        4.6 &       0.64 &       0.29 &          8 &          5 \\ \bottomrule
		\end{tabular}  

\label{tab:technicalchp}
\end{table*}
		\begin{table*}
		\caption{Technical parameters of the auxiliary unit and storages}

		\footnotesize
		\centering
		\begin{tabular}{p{0.1pt}p{15pt}p{22pt}p{8pt}p{8pt}p{8pt}p{8pt}p{8pt}p{8pt}p{8pt}lll}\toprule
		
		\textbf{} & \multicolumn{ 2}{c}{{Aux. boiler}} &              \multicolumn{ 4}{c}{{Thermal storage}} &                           \multicolumn{ 6}{c}{{Biomass storage}} \\ \cmidrule(r){2-3} \cmidrule(r){4-7} \cmidrule(r){8-13}
		
		&    $E^{\text{AUX}}$               &    $\overline{Q^{\text{AUX}}}$    &    $S^{0}$  &    $S^{F}$  &    $\overline{S}  $  &    $\underline{S}$  &        $\Delta^{\text{F}}$   &      $\Delta^{0}$        &       $\Delta^\text{W}$ & $\overline{\Delta} $  &    $\underline{\Delta}_t$  & $E^{\text{B}}$ \\\midrule
		
		A &       0.97 &         15 &          5 &          3 &          7 &          0 &           35 &        500 &         24 &  20000   &   4000 (20-45) & 4.9971\\
		&&&&&&&& &&&  2000 &\\
		B &       0.99 &         30 &        6.5 &        4.5 &        9.5 &          0 &           70 &        850 &         24 & 35000 &   7000 (20-45) & 4.9971\\
		&&&&&&&& &&&  3500 & \\\bottomrule
		
	\end{tabular}  
	
	\label{tab:technicalaux}
\end{table*}

\begin{table*}
	\centering
	\footnotesize
	\caption{Cost parameters}
	
	\begin{tabular}{rrrrrrrrrr}
		\toprule
	 &                                 \multicolumn{ 5}{c}{{CHP}} &      \multicolumn{ 3}{c}{{Aux. boiler}} & \multicolumn{ 1}{c}{{Storage}} \\\cmidrule(r){2-6} \cmidrule(r){7-9}
		
		&    $C^{\text{CHP}}$    &    $C^{\text{SU}}$            &    $C^{\text{SD}}$            &    $T^{\text{EP}}$            &    $I$                        &    $C^{\text{O\&M}}_{Aux}$    &    $T^{\text{AUX}}$           &    $T^{\text{CO}_2}$     & $C^\text{I}$     \\ \midrule
		
		A &      19.85 &      14250 &          0 &      55.62 &      20.25 &       0.07 &      28.22 &       6.34 &    0.0002 \\
		
		B &      20.32 &      16870 &          0 &      55.62 &      20.25 &        0.5 &      52.07 &          0 &      0.0002 \\
		\bottomrule
	\end{tabular}  
	\label{tab:cost}
\end{table*}

\begin{table*}
	\centering
	\footnotesize
	\caption{Biomass contract data}
	\begin{tabular}{crrrrrrrrrr}
		\toprule
		Contract	&	$C_{j}^{\text{B}}$	&$C_{j}^{\text{B}+}$	& $C_{j}^{\text{B}-}$	&	 $O_{j}^{+}$	&	$O_{j}^{-}$	&$\overline{B_{j}}$ &$\underline{B_{j}}$	&$F_{j}$	&	$\overline{N_{j}}$ 	&$\underline{N_{j}}$ \\ \hline
		1	&	150.8	&	0	&	0	&	0	&	0	&	19000	&	18000	&	2016	&	4	&	4 \\
		2	&	156.4	&	0	&	0	&	0	&	0	&	17000	&	12000	&	1344	&	5	&	2 \\
		3	&	170.83	&	0	&	0	&	0	&	0	&	15000	&	11000	&	1008	&	8	&	4 \\
		4	&	181.31	&	30.56	&	30.56	&	0.1	&	0.1	&	12000	&	8000	&	504	&	17	&	15 \\
		5	&	181.43	&	24.45	&	24.45	&	0.15	&	0.15	&	12000	&	8000	&	504	&	15	&	15 \\
		6	&	183.59	&	30.56	&	30.56	&	0.25	&	0.25	&	5100	&	2380	&	336	&	25	&	24 \\ 
		7	&	183.43	&	36.67	&	36.67	&	0.25	&	0.25	&	5100	&	2380	&	336	&	25	&	15 \\
		8	&	201.89	&	18.34	&	18.34	&	0.5	&	0.5	&	1200	&	1200	&	168	&	50	&	50 \\
		9	&	202.17	&	18.34	&	18.34	&	0.5	&	0.5	&	1200	&	1000	&	168	&	50	&	25 \\
		10	&	204.29	&	28.12	&	28.12	&	0.5	&	0.5	&	850	&	850	&	120	&	60	&	50 \\
		11	&	202.24	&	28.12	&	28.12	&	0.65	&	0.65	&	850	&	500	&	120	&	60	&	30 \\
		12	&	202.05	&	12.22	&	12.22	&	0.75	&	0.75	&	350	&	100	&	48	&	100	&	80 \\
		13	&	202.64	&	12.22	&	12.22	&	0.75	&	0.75	&	350	&	100	&	48	&	100	&	50 \\
		\bottomrule
	\end{tabular}  
	\label{tab:contracts}
\end{table*}

\subsection{Scenario generation}
Apart from the deterministic parameters mentioned in the previous section, we have to handle uncertainty regarding heat demands, gas prices and electricity prices to be used in the optimization. Since both municipalities are within the same bidding region in Nordpool (DK1) and the same gas trading region, the electricity and natural gas prices are identical. However, differences exist regarding the heat demand. We use historical data from 1st June 2011 to 31st May 2016  for electricity prices, natural gas prices and heat demands. Based on this data, different techniques for scenario generation are implemented.  
The resulting scenarios and expected values depend on the municipality due to the different auxiliary boilers and heat consumption in previous years. Furthermore, the input time series varies with the phase of the solution approach regarding time scales and need for scenarios. The scenario generation for both phases is described in Appendix \ref{sec:appendix_a}.

\subsection{Evaluation of solution approach} 

To evaluate our solution approach, we have to obtain the costs under different realizations of the uncertainty. We use 11 samples, i.e, 11 different realizations of uncertainty, for each municipality. Sample 0 is the actual realization of the heat demand, electricity prices and gas prices from 1st June 2016 to 31st May 2017. The remaining 10 samples (from 1 to 11) are a composite of different real data sets obtained from the same sources as the previous data. The electricity and gas prices are obtained from real data of 2015, 2016 and 2017 from other regions in Nordpool and other European hubs, respectively. The heat consumption is obtained from other municipalities in the Aarhus district heating system and scaled to the size of the system capacity accordingly.

The scenarios used in the evaluations are based on a combination of past data and time series forecasts (see description in Appendix \ref{sec:appendix_a}). To decide the biomass supply contracts, we use expected values for electricity and gas prices based on the last five years, while the heat demand is modeled using five scenarios resembling the heat demand from the last five years. In the operational planning problem, the scenarios for electricity prices and heat demand consist of a time series forecast for the next week while the corresponding historical values of electricity prices and heat demand in previous years are taken for the remaining weeks of the receding horizon (denoted as method \emph{F1}). See Appendix \ref{sec:appendix_b} for a comparison of different scenario generation methods.

Evaluating one sample with a configuration of our method requires to extend Algorithm \ref{alg:approach} by one step. Each week after the operational problem is solved (line 5 in Algorithm \ref{alg:approach}), we fix the first-stage decisions and solve the model  using the realizations of the uncertainty of the first week. Thus, we obtain the real costs for the first week and the initial status for the next week.

\section{Experimental results}\label{sec:results}
For the experimental evaluation, we implemented Algorithm \ref{alg:approach} using Python 3.5.1 and Gurobi 7.0.1 (default parameters). All experiments are run on Intel Xeon Processor X5550 with 24 GB RAM. The objective values in this section comprise the real costs summed over all weeks in the year.

\subsection{Analysis of receding horizon length}\label{sec:results_horizon}

\begin{table*}
	\footnotesize
	\centering
	\caption{Objective value and penalty costs [x100,000\euro] for different lengths of receding horizon}
	
	\begin{tabular}{rrrrrrrrrr}

		\toprule
		\multicolumn{2}{c}{Sample}     &            \multicolumn{ 4}{c}{Objective} &          \multicolumn{ 4}{c}{Penalty $q^{\text{miss}}$} \\  \cmidrule(r){1-2} \cmidrule(r){3-6} \cmidrule(l){7-10}
		
		&    &  1 &  2 &  3 &  4 &1 &  2 &  3 &  4 \\\midrule
		
		\multirow{11}{*}{\rotatebox[origin=l]{90}{Municipality A}}  &          0 &     91.635 & {\bf 84.505} &     84.548 &     84.687 &      7.576 &      1.376 &      1.376 &      1.376 \\
		
		&          1 &     96.483 &     84.571 &     84.540 & {\bf 84.505} &     10.993 &      0.000 &      0.000 &      0.000 \\
		
		&          2 &     87.707 &     81.793 &     81.810 & {\bf 81.710} &      4.653 &      0.000 &      0.000 &      0.000 \\
		
		&          3 &     96.915 & {\bf 84.364} &     84.407 &     84.389 &     11.521 &      0.000 &      0.000 &      0.000 \\
		
		&          4 &     90.140 & {\bf 82.205} &     82.455 &     82.394 &      6.582 &      0.000 &      0.000 &      0.000 \\
		
		&          5 &     86.659 &     83.756 &     83.735 & {\bf 83.707} &      2.130 &      1.010 &      1.010 &      1.010 \\
		
		&          6 &    104.655 &     82.764 &     82.660 & {\bf 82.511} &     21.334 &      0.000 &      0.000 &      0.000 \\
		
		&          7 &     89.979 & {\bf 87.009} &     87.027 &     87.031 &      5.832 &      2.904 &      2.904 &      2.904 \\
		
		&          8 &     89.700 &     82.035 & {\bf 81.846} &     81.847 &      7.187 &      0.000 &      0.000 &      0.000 \\
		
		&          9 &     84.836 & {\bf 82.469} &     82.588 &     82.609 &      1.192 &      0.000 &      0.000 &      0.000 \\
		
		&         10 &     84.900 & {\bf 83.515} &     83.635 &     83.580 &      0.439 &      0.131 &      0.131 &      0.131 \\\hline
		
		\multirow{11}{*}{\rotatebox[origin=l]{90}{Municipality B}}  &          0 &    400.734 & {\bf 169.017} &    169.186 &    169.682 &    231.804 &      1.297 &      1.297 &      1.297 \\
		
		&          1 &    246.238 & {\bf 174.582} &    175.025 &    175.056 &     71.916 &      0.000 &      0.000 &      0.000 \\
		
		&          2 &    267.247 &    167.181 & {\bf 167.022} &    167.371 &     95.517 &      0.000 &      0.0 00 &      0.000 \\
		
		&          3 &    286.542 &    175.332 &    175.059 & {\bf 175.059} &    111.651 &      0.000 &      0.000 &      0.000 \\
		
		&          4 &    318.557 &    169.592 &    169.334 & {\bf 168.766} &    147.841 &      0.000 &      0.000 &      0.000 \\
		
		&          5 &    386.050 &    170.939 &    170.776 & {\bf 170.597} &    216.553 &      1.335 &      1.335 &      1.335 \\
		
		&          6 &    171.352 &    170.397 & {\bf 170.390} &    171.026 &      0.000 &      0.000 &      0.000 &      0.000 \\
		
		&          7 &    275.865 &    172.744 &    172.742 & {\bf 172.672} &    104.186 &      0.000 &      0.000 &      0.000 \\
		
		&          8 &    169.469 &    168.003 &    167.982 & {\bf 167.525} &      0.000 &      0.000 &      0.000 &      0.000 \\
		
		&          9 &    323.960 & {\bf 171.312} &    171.693 &    171.845 &    152.953 &      0.000 &      0.000 &      0.000 \\
		
		&         10 &    397.496 &    169.546 &    169.399 & {\bf 168.974} &    225.878 &      0.000 &      0.000 &      0.000 \\
		\bottomrule
	\end{tabular}  
	\label{table:results_horizon}
	
\end{table*}

Table \ref{table:results_horizon} shows the objective value and penalty costs using scenarios generated by method \emph{F1} for different lengths of the receding horizon, namely one, two, three and four weeks. Note that in no case, penalty costs for exceeding the biomass storage capacity occurred and therefore those are omitted from the table.  The most important result is that the objective values drastically improves by including at least a second week into the optimization. For the most part, this is due to the reduction in penalty costs for not fulfilling the heat demand (see Table \ref{table:results_horizon}). This can be explained by the opportunity of using of options. Indeed, if the receding horizon already considers scenarios for weeks apart from the current week, we make use of this information now. If the biomass contract selection (phase 1) scheduled a delivery only in the current week, but not in the next week, having a longer planning horizon can be beneficial. If we only consider the current week, we may not make use of an upward option, because it is not needed now. However, if the scenarios for the next week(s) show a trend with a higher heat demand than expected, we can get more biomass than scheduled now instead of running out of storage and missing the demand.

For a receding horizon of more than one week, the results are quite similar. The maximum deviation between costs for the different lengths is 82548 for municipality B in sample 4, which is in total approximately 0.48\% higher costs. For all other cases, the relative and absolute difference is less. In some cases, a longer horizon can lead to slightly poorer results due to the fact, that the heat demand is still uncertain and we may make use of an upward or downward option that corrects the delivery amount according to the uncertain scenarios. If the scenarios show a wrong trend in later weeks, it can be more beneficial to just include a second week (e.g. sample 0, mun. A). 
The penalty cost for missing the heat demand is $\phi^{\text{Miss}}=10000$ [\euro/MWh], which means we miss at most 29.0356 MWh of heat in sample 7 for municipality A in a whole year. In all cases with penalty cost, the missing demand occurs in periods with an exceptionally high demand close to the capacity of the system. Those very high demands are often not covered by the scenarios and therefore wrong planning decisions may cause a shortage of biomass and a penalization for not satisfying the heat demand.  Note that in practice a lack of supply in the district heating network would never occur, because the heat producer can gradually decrease the supply temperature or reduce the water flow to increase the demand covered. However, these cases must be avoid and therefore we penalize them  in the objective .

\subsection{Stochastic programming vs. expected value solution}\label{sec:results_stovsdet}

\begin{table*}[h]
	\footnotesize
	\centering
	\caption{Comparison stochastic programming (Sto.) vs. expected value solution (Exp.) [x100,000\euro]. The maximum, minimum and average values are based on 2 to 4 weeks horizon}
	\begin{tabular}{rrrrrrrr}
		\toprule
		\multicolumn{2}{c}{Sample} &    Max. Sto. &    Min. Exp. &      Delta &    Avg. Sto. &    Avg. Exp. &      Delta \\\midrule
		
		\multirow{12}{*}{\rotatebox[origin=l]{90}{Municipality A}}&          0 & {\bf 84.687} &     84.975 &     0.34\% & {\bf 84.580} &     85.117 &     0.63\% \\
		
		&          1 & {\bf 84.571} &     85.331 &     0.89\% & {\bf 84.539} &     85.343 &     0.94\% \\
		
		&          2 & {\bf 81.810} &     82.010 &     0.24\% & {\bf 81.771} &     82.208 &     0.53\% \\
		
		&          3 & {\bf 84.407} &     85.205 &     0.94\% & {\bf 84.387} &     85.315 &     1.09\% \\
		
		&          4 & {\bf 82.455} &     82.853 &     0.48\% & {\bf 82.352} &     83.008 &     0.79\% \\
		
		&          5 & {\bf 83.756} &     84.354 &     0.71\% & {\bf 83.733} &     84.462 &     0.86\% \\
		
		&          6 & {\bf 82.764} &     83.210 &     0.54\% & {\bf 82.645} &     83.253 &     0.73\% \\
		
		&          7 & {\bf 87.031} &     87.744 &     0.81\% & {\bf 87.022} &     87.769 &     0.85\% \\
		
		&          8 & {\bf 82.035} &     82.445 &     0.50\% & {\bf 81.909} &     82.501 &     0.72\% \\
		
		&          9 & {\bf 82.609} &     83.381 &     0.93\% & {\bf 82.556} &     83.475 &     1.10\% \\
		
		&         10 & {\bf 83.635} &     84.006 &     0.44\% & {\bf 83.577} &     84.072 &     0.59\% \\
		
		&       Avg. & {\bf 83.615} & { 84.138} & { 0.62\%} & {\bf 83.552} & { 84.229} & { 0.80\%} \\
		
		\hline
		\multirow{12}{*}{\rotatebox[origin=l]{90}{Municipality B}} &          0 & {\bf 169.682} &    172.264 &     1.50\% & {\bf 169.295} &    172.264 &     1.72\% \\
		
		&          1 & {\bf 175.056} &    177.669 &     1.47\% & {\bf 174.888} &    177.694 &     1.58\% \\
		
		&          2 & {\bf 167.371} &    170.974 &     2.11\% & {\bf 167.191} &    170.986 &     2.22\% \\
		
		&          3 & {\bf 175.332} &    176.147 &     0.46\% & {\bf 175.150} &    176.183 &     0.59\% \\
		
		&          4 & {\bf 169.592} &    171.702 &     1.23\% & {\bf 169.231} &    171.720 &     1.45\% \\
		
		&          5 & {\bf 170.939} &    172.536 &     0.93\% & {\bf 170.770} &    172.546 &     1.03\% \\
		
		&          6 & {\bf 171.026} &    171.303 &     0.16\% & {\bf 170.604} &    171.328 &     0.42\% \\
		
		&          7 & {\bf 172.744} &    174.844 &     1.20\% & {\bf 172.719} &    174.844 &     1.22\% \\
		
		&          8 & {\bf 168.003} &    168.231 &     0.14\% & {\bf 167.837} &    168.234 &     0.24\% \\
		
		&          9 & {\bf 171.845} &    172.635 &     0.46\% & {\bf 171.617} &    172.641 &     0.59\% \\
		
		&         10 &    169.546 & {\bf 169.456} &    -0.05\% & {\bf 169.306} &    169.703 &     0.23\% \\
		
		&       Avg. & {\bf 171.012} & { 172.524} & { 0.87\%} & {\bf 170.783} & { 172.558} & { 1.03\%} \\
		\bottomrule
	\end{tabular}  
	
	\label{tab:results_sto_vs_det}
\end{table*}

To show the benefit of using stochastic programming instead of using an expected value approach, we compare the results in Table \ref{tab:results_sto_vs_det}. We limit our results to \emph{F1} scenarios and 2 to 4 weeks of receding horizon. The first three columns in Table \ref{tab:results_sto_vs_det} compare the worst case among the three lengths of receding horizon in the stochastic approach with the best case of the expected value solution. The worst case stochastic solution gives on average a 0.75\% better solution and it dominates in all cases except one (municipality B, sample 10). When comparing the average objective values, the stochastic approach improved the results on average by 0.80\% for municipality A and 1.03\% for municipality B. Although the improvement is not relatively large, in absolute terms it results in saving on average \euro{67743} and \euro{177583}, respectively.

The improvement can be explained by the fact that the stochastic solution makes use of options while the expected value solution does not contract any options (see Section \ref{sec:furtherresults} for more details).

\subsection{Interpretation of results for real data from 2016-2017}\label{sec:furtherresults}

In this section, we describe the results of the contract selection and operational planning in more detail. As an example, we analyze sample 0 for both municipalities, which contains the real data from 1st June 2016 to 31 May 2017. We would like to point out that the conclusions drawn in this section coincides with the observations from the other samples.

\begin{figure*}
	\begin{subfigure}[t]{.5\textwidth}
		\centering
		\includegraphics[width=\textwidth]{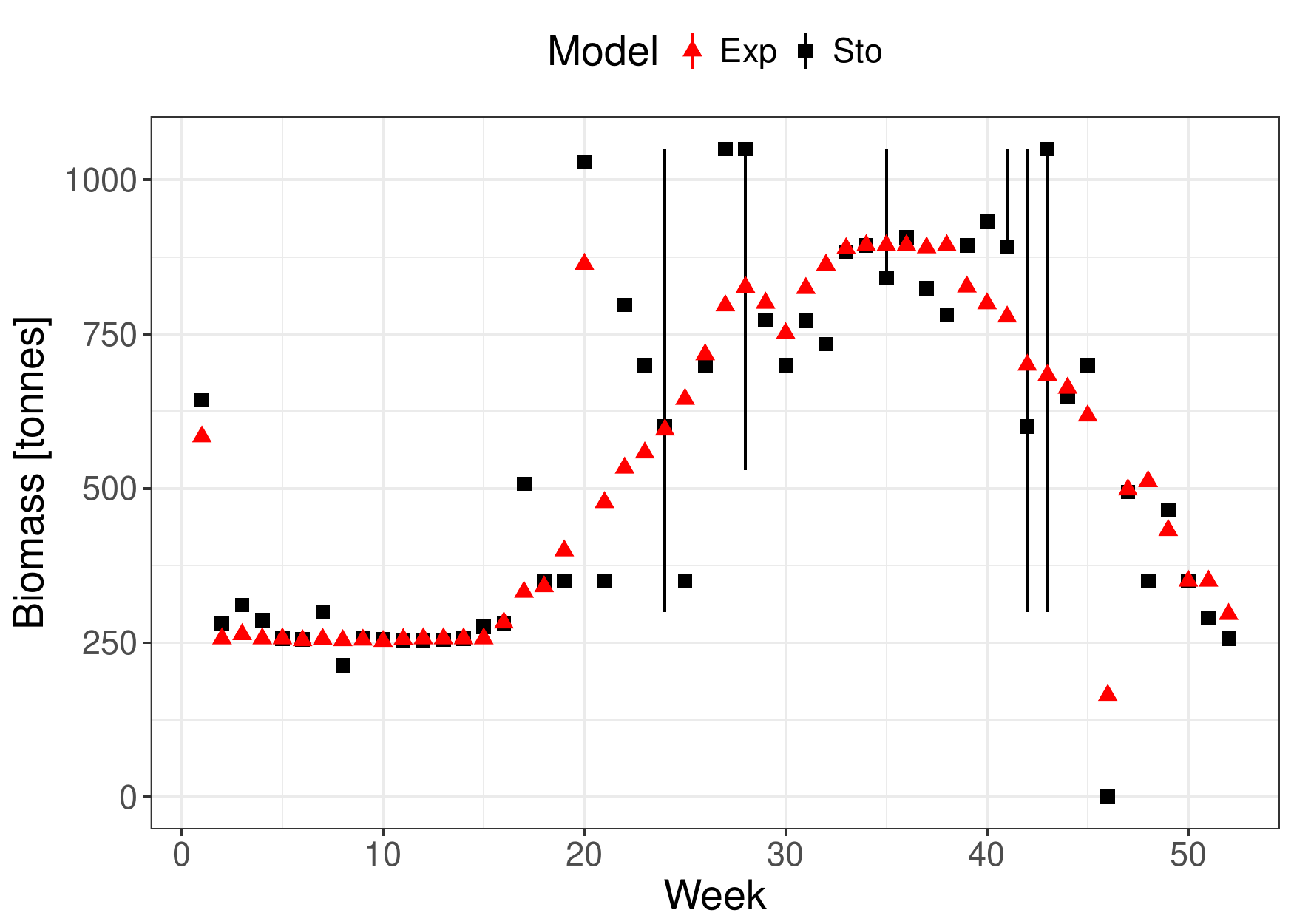}
		\caption{Municipality A (selected contract: 12)}
	\end{subfigure}
	\begin{subfigure}[t]{.5\textwidth}
		\centering
		\includegraphics[width=\textwidth]{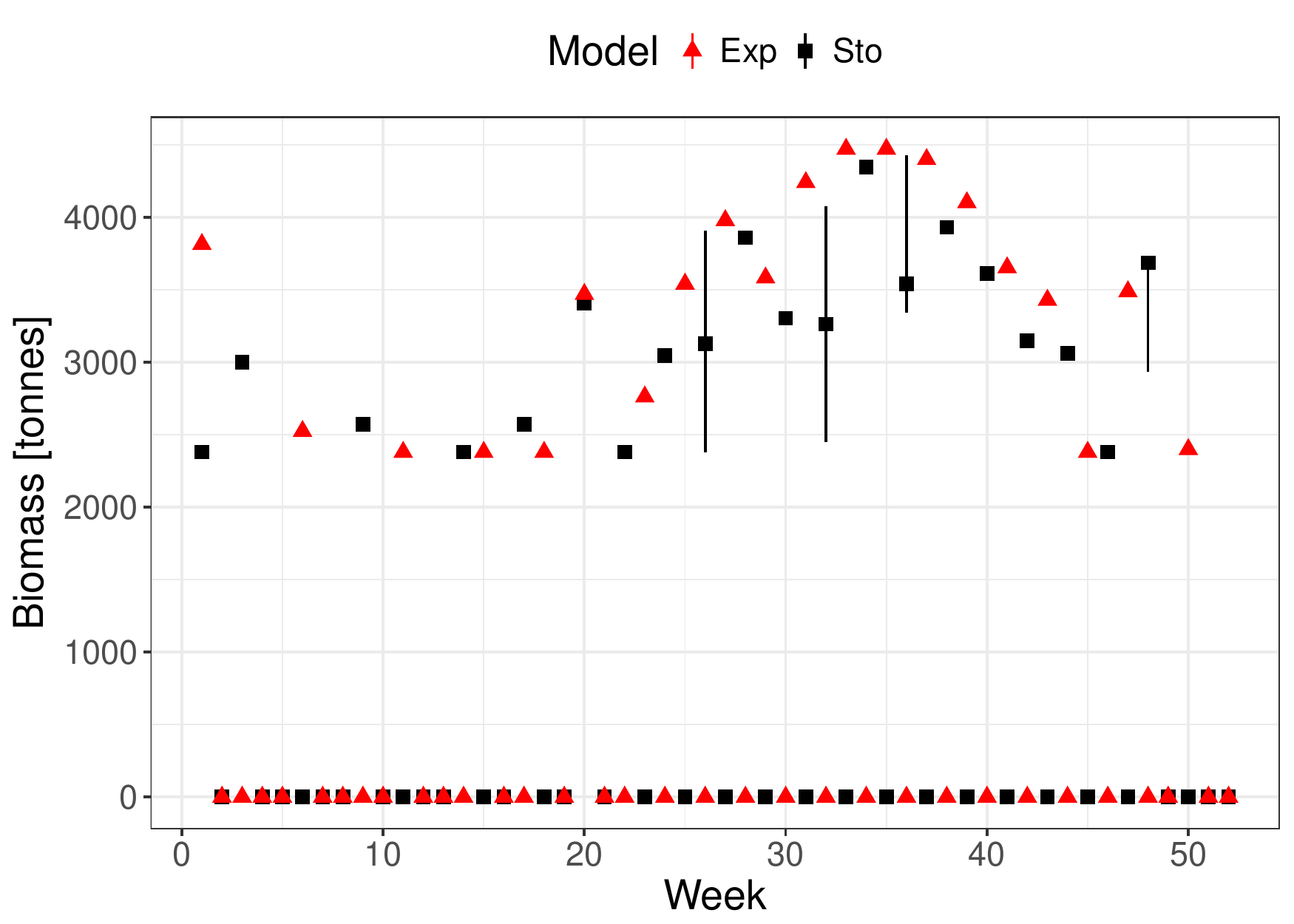}
		\caption{Municipality B (selected contract: 7)}
	\end{subfigure}
	\caption{Biomass contracts from 1 June 2016 to 31 May 2017 }
	\label{fig:schedule}
\end{figure*}

Figure \ref{fig:schedule} shows the selected biomass contracts for municipalities A and B in the stochastic (\textit{Sto}) and expected value solution (\textit{Exp}), respectively. The contracts are valid for all samples. The points show the contracted biomass amount and the vertical lines that extent from some of the crosses are the amount of upward and downward options bought.  
We see that only the solutions obtained by the stochastic approach make use of options. As the deterministic solution has no scenarios and assumes the expected values of uncertain parameters as deterministic, the contracts are selected in such a way that the solution fits these expected values. Thus, no use of options is reasonable in this case. However, when other biomass amounts are needed in the course of the year, the options contracted in the stochastic solution bring an advantage and reduce the overall cost (see Table \ref{tab:results_sto_vs_det}). From Figure \ref{fig:schedule}  also the difference in the delivery patterns for the two municipalities can be seen. The selected contract for municipality A (contract 12) has smaller amounts but more frequent deliveries. Whereas the selected contract for municipality B has larger amounts and less deliveries, which relates also to the higher heat demand in municipality B.

\begin{figure*}
	\begin{subfigure}[t]{.5\textwidth}
		\centering
		\includegraphics[width=\textwidth]{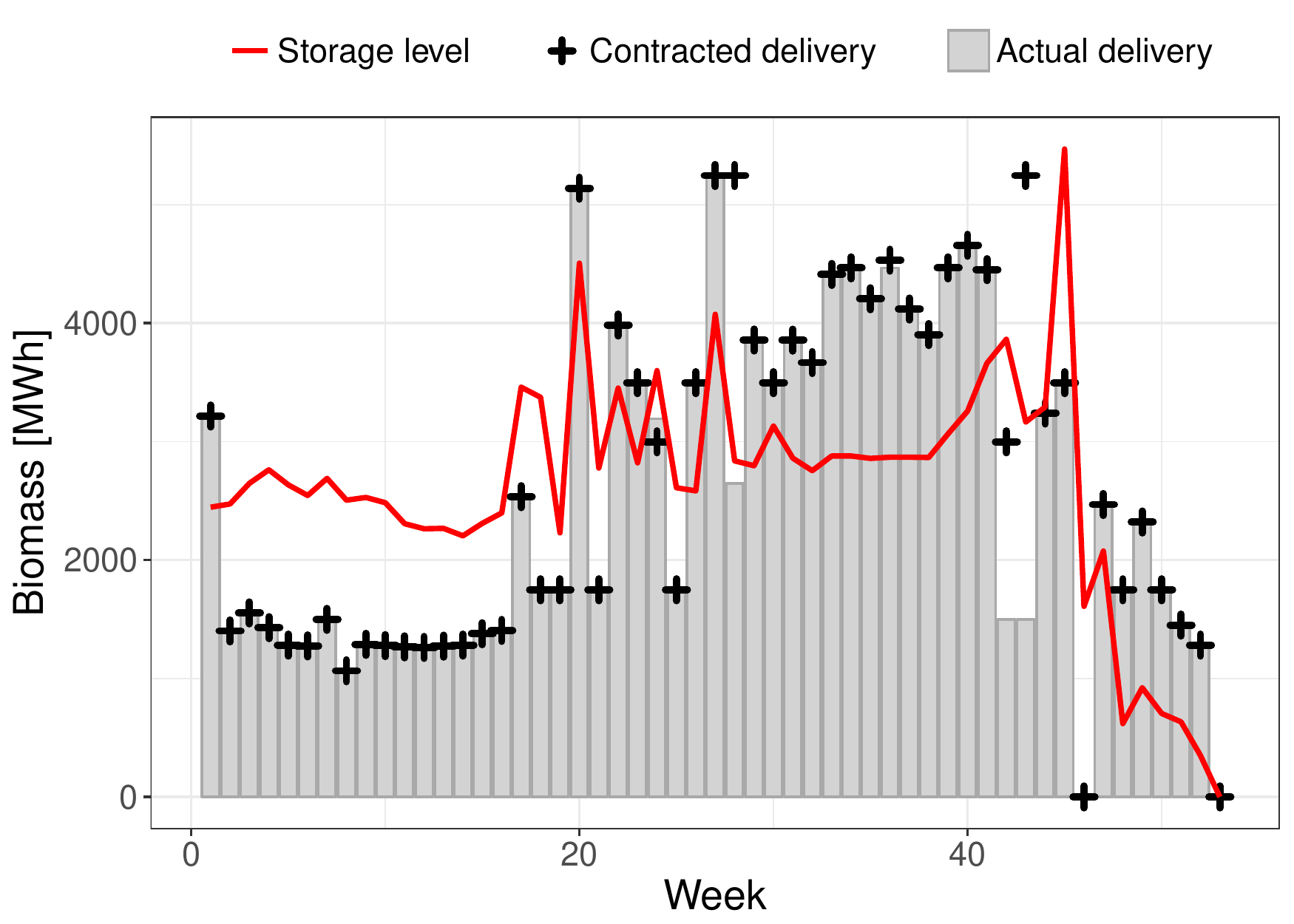}
		\caption{Municipality A}
	\end{subfigure}
	\begin{subfigure}[t]{.5\textwidth}
		\centering
		\includegraphics[width=\textwidth]{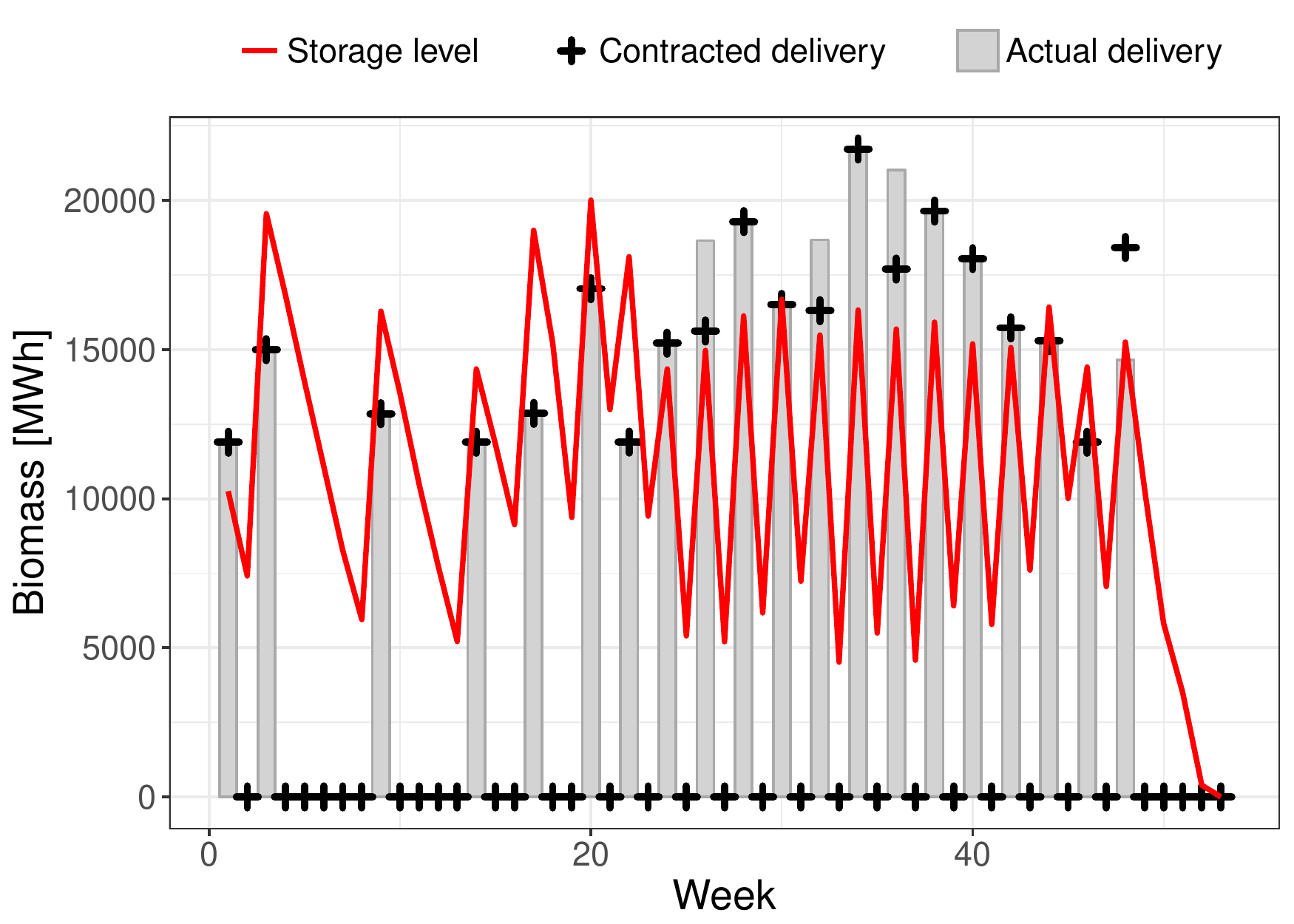}
		\caption{Municipality B}
	\end{subfigure}
	\caption{Biomass storage level and deliveries for the real realization of uncertainties from 1st June 2016 to 31 May 2017 (based on \emph{F1} scenarios and two weeks of receding horizon)}
	\label{fig:biomass_storage}
\end{figure*}

The actual delivery amounts, i.e., after making use of options, and the biomass storage level are depicted in Figure \ref{fig:biomass_storage} for the real data of the year 2016 to 2017. The amounts are cumulated per week. Furthermore, the contracted delivery amount is depicted to show if the options are actually used in the course of the year. For both municipalities the operational problem uses both upward and downward options, for example, week 24 (downward) and 42 (upward) in municipality A or week 26 (downward) and 48 (upward) in municipality B.

\begin{figure*}
	\begin{subfigure}[t]{.5\textwidth}
		\centering
		\includegraphics[width=\textwidth]{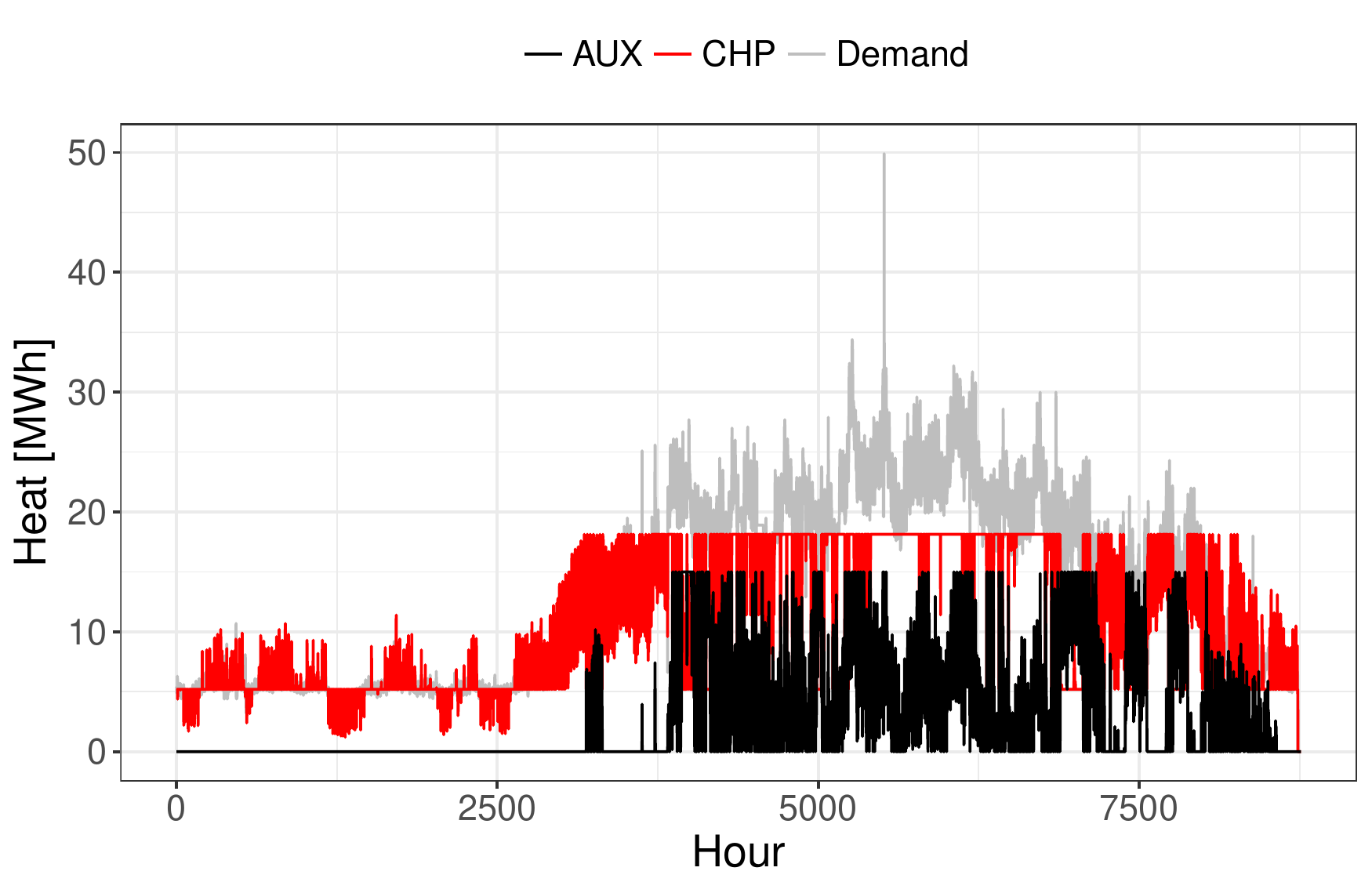}
		\caption{Municipality A}
		\label{fig:heat_real_a}
	\end{subfigure}
	\begin{subfigure}[t]{.5\textwidth}
		\centering
		\includegraphics[width=\textwidth]{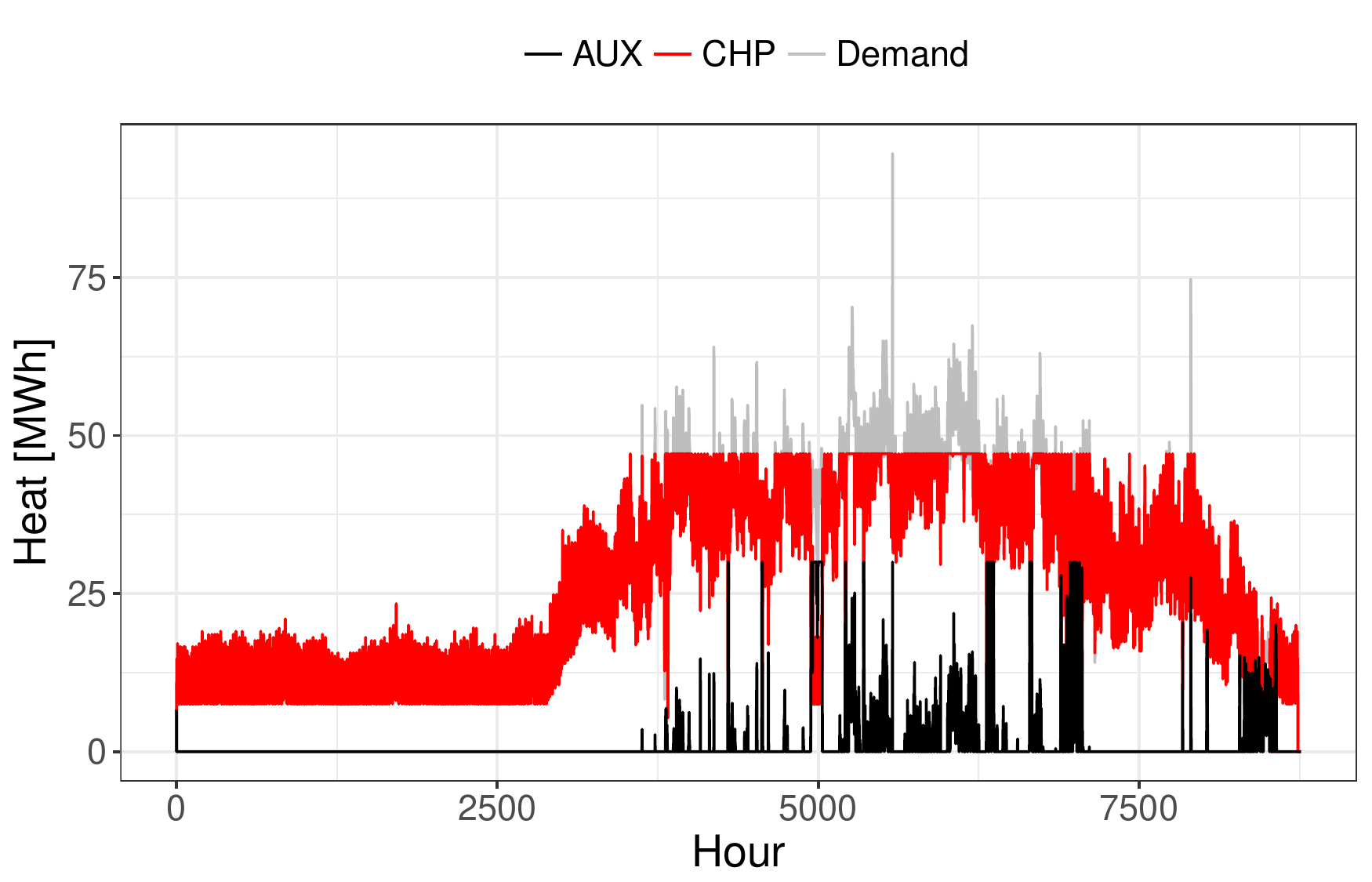}
		\caption{Municipality B}
		\label{fig:heat_real_b}
	\end{subfigure}
	\caption{Heat production for real realization of uncertainties from 1st June 2016 to 31 May 2017 (based on \emph{F1} scenarios and two weeks of receding horizon)}
	\label{fig:heat_real}
\end{figure*}

The heat production from June 2016 to May 2017 for municipality A and B is shown in Figure \ref{fig:heat_real}. In both cases the heat demand was always fulfilled and the production  follows similar behavior. At start of the season, the demand can be covered by the biomass-fired CHP. During the winter periods with a high demand, the gas boiler is used in addition to the CHP to cover the heat demand. Furthermore, at the end of the season the boiler is used more often as in the beginning of the season due to a slightly higher demand and the biomass contract decisions contracting less biomass in the end of the season. 

\subsection{Runtime analysis}
\begin{figure*}
	\begin{subfigure}{.5\textwidth}
		\centering
		\includegraphics[width=1\textwidth]{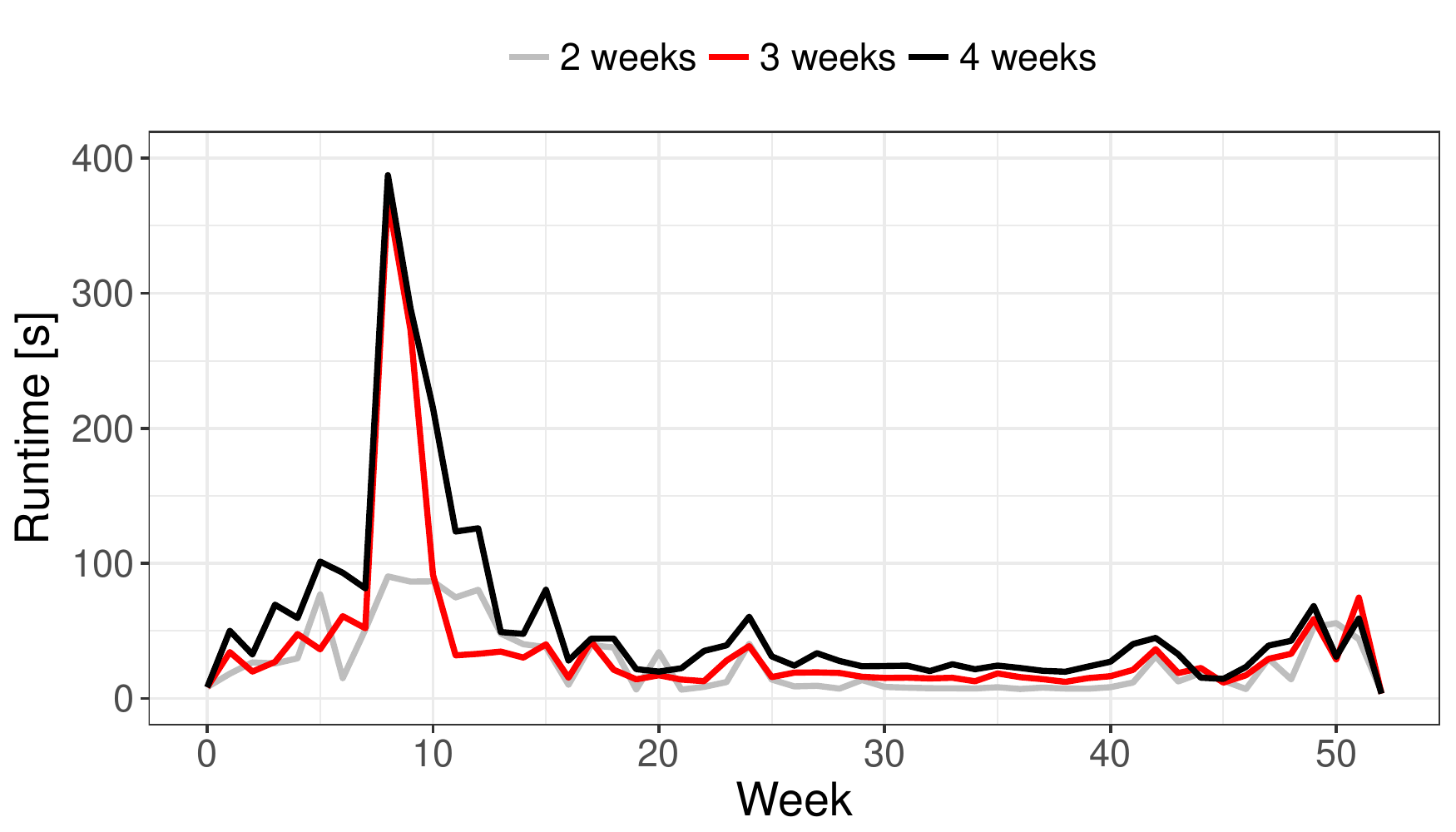}
		\caption{Municipality A } \label{fig:runtime_a}
	\end{subfigure}
	\begin{subfigure}{.5\textwidth}
		\centering
		\includegraphics[width=1\textwidth]{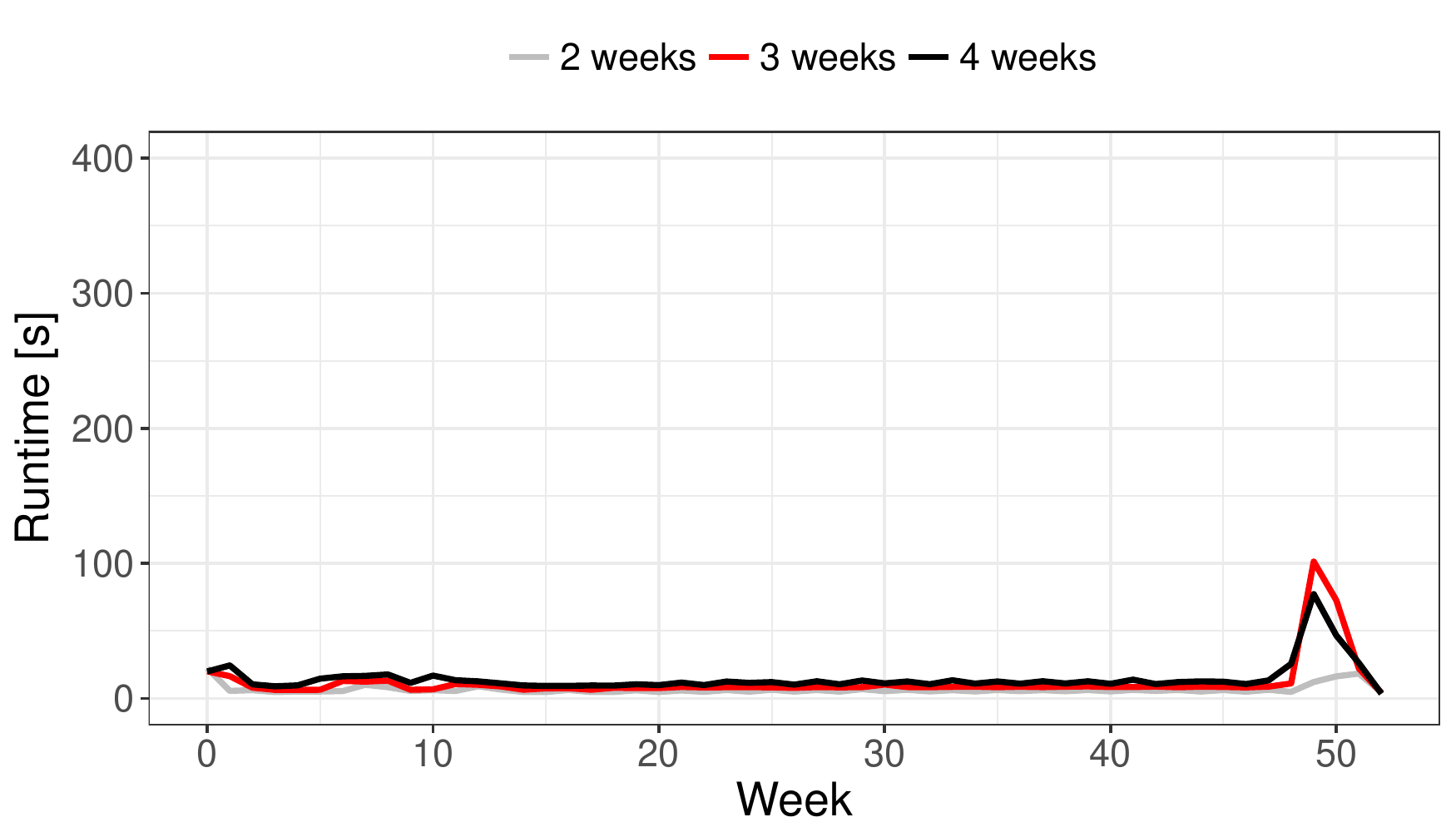}
		\caption{Municipality B } \label{fig:runtime_b}
	\end{subfigure}
	\caption{Average runtimes per week (week 0 corresponds to biomass contract selection)}
	\label{fig:runtimes}
\end{figure*}

Figure \ref{fig:runtimes} shows the runtimes for different lengths of the receding horizon averaged over the 11 samples. The corresponding MIP model sizes for the biomass selection and the operational planning problem with different lengths of the receding horizon are given in \ref{tab:modelsizes}. Note that the model size for each week in the operational planning phase is the same throughout the planning horizon. Therefore, the model sizes depend only on the length of the receding horizon. For most of the cases, the runtime to solve the operational model for one week is less than 60 seconds. Also, the biomass contract selection model is solved in less than 20 seconds for both municipalities (see week 0 in Fig. \ref{fig:runtimes}). The runtime slightly increases with a longer receding horizon, but not significantly.

\begin{table}
	\caption{Model sizes}
	\footnotesize
	\begin{tabular}{lrrrrr}\toprule
		& Cont. var. & \multicolumn{2}{c}{Int. var.  (thereof bin. var.)}& Constraints & NZs\\
		\midrule
		Biomass selection & 10,878 & 689 & (13) & 1,3812 & 40,479\\
		Operational - 1 week & 15,125 & 3,360 & (3,360) & 39,495 & 152,180\\
		Operational -	2 weeks & 30,250 & 6,720 & (6,720) & 65,530 & 284,900 \\
		Operational - 3 weeks & 45,375 & 10,080 & (10,080) & 91,615 & 417,620\\
		Operational - 4 weeks & 60,500 & 13,440 & (13,440) & 117,675 & 550,340 \\\bottomrule
	\end{tabular}
	\label{tab:modelsizes}
\end{table}

For the few cases with a high runtime the average lies below 400 seconds (see \ref{fig:runtime_a}), which is short enough for a weekly planning problem to be used in practice. The weeks with higher runtime relate to samples where the heat demand is higher than expected in the biomass contract selection phase, which leads to a shortage of biomass in the subsequent weeks (in the beginning of the year in municipality A and in the end of the year in municipality B). Due to this shortage the model tries to avoid penalties for getting below the safety storage level while producing as much as possible with the CHP to get income from the electricity market. As the production is not possible in all hours, the model has to select the hours with highest expected electricity prices making it harder for the solver to find the best solution as the electricity prices are close to each other.

\section{Summary and outlook}
\label{conclusion}

In this work, we propose a solution approach that optimizes the biomass supply planning for a large-scale CHP producer using biomass. The decision-making process is divided into two phases both using two-stage stochastic programs. The first model, named \textit{biomass contract selection}, is solved for a long-term horizon with weekly periods and configures the contracts from a set of biomass suppliers. Those decisions are used in the second model, named \textit{operational planning}, to optimize the heat production. This solutions approach corresponds to the planning process in practice.
We evaluate our method on two case studies with realistic requirements and historical data to create scenarios. We analyze several scenario generation possibilities to create the scenarios based on past data and different forecasting tools. Our analysis investigates the results obtained for 11 samples of realizations of uncertainty. 

The results reflect that the use of a receding horizon improves the solution obtained due to a better operation of both heat and biomass storages. However, as a result of the forecast uncertainty, very long receding horizons may not improve the results. Furthermore, we show that applying stochastic programming is required to make use of the options, yielding better results than in the expected value case where no options are purchased. 

We envision four future research directions. First, further uncertainties regarding the delivery of biomass such as amount and quality variations could be included in a supply chain planning model. Second, the configuration of our algorithm can be investigated further to determine the length of the receding horizon in a better way and improve the results. Third, an economic analysis of the options can be made to assess their benefit for the entire supply chain. That is from both supplier's and producer's points of view. Finally, the comparison of different long-term forecasting tools with the use of data from previous years to create long-term scenarios is another future research direction. 

\section*{Acknowledgements}
The work of Ignacio Blanco, Daniela Guericke and Henrik Madsen is funded by Innovation Fund Denmark through the CITIES research center (no. 1035-00027B). The work of Juan M. Morales is partly funded by the Spanish Research Agency through project
ENE2017-83775-P (AEI/FEDER, UE), and by the Research Funding Program for Young Talented
Researchers of the University of Malaga through project PPIT-UMA-B1-2017/18. The authors would like to thank Ørsted A/S for their valuable input
and comments to this work.


\begin{appendices}

\section{Scenario generation}\label{sec:appendix_a}
In this section, we describe the different approaches used for scenario generation in biomass contract selection and operational planning problem, respectively.

\subsection{Biomass contract selection} In phase one of the solution approach, scenarios for the heat demand and the expected value for auxiliary boiler costs and electricity prices are part of the model. In this tactical planning problem, we use the heat consumption of the five previous years (i.e. 1st June 2011 - 31st May 2016) from summer to summer of the respective community as heat demand scenarios ($D_{t,\omega}$) resulting in five scenarios. The probability for each scenario is determined based on the year while giving a higher probability to more recent years (first three years: 0.15, last two years: 0.275).  

\begin{table*}[t]
	\footnotesize
	\centering
	\caption{Minimum objective value [x100,000\euro] for each mode of scenario generation. The minimum refers to the lowest objective value of 1, 2, 3 or 4 weeks of receding horizon }
	
	\begin{tabular}{p{0.009\textwidth}p{0.06\textwidth}p{0.06\textwidth}p{0.06\textwidth}p{0.06\textwidth}p{0.06\textwidth}p{0.06\textwidth}p{0.06\textwidth}p{0.06\textwidth}p{0.06\textwidth}p{0.06\textwidth}}
		\toprule
		{\bf S.} &                      \multicolumn{ 5}{c}{{\bf Municipality A}} &                      \multicolumn{ 5}{c}{{\bf Municipality B}} \\
		\cmidrule(r){2-6} \cmidrule(l){7-11}
		
		&    { P} &   { F1} & { F1+P} &   { F2} & { F2+P} &    { P} &   { F1} & { F1+P} &   { F2} & { F2+P} \\
		\midrule
		0 &      84.61 & {\bf 84.51} &      84.63 &      84.88 &      85.08 &     169.25 & {\bf 169.02} &     169.20 &     169.28 &     169.27 \\
		
		1 &      84.62 & {\bf 84.50} &      84.61 &      84.78 &      84.59 &     174.94 & {\bf 174.58} &     174.66 &     174.99 &     174.63 \\
		
		2 &      81.69 &      81.71 & {\bf 81.67} &      81.72 &      81.79 &     167.19 & {\bf 167.02} &     167.30 &     167.69 &     167.33 \\
		
		3 &      84.55 & {\bf 84.36} &      84.46 &      84.53 &      84.48 &     175.23 & {\bf 175.06} &     175.20 &     175.23 &     175.18 \\
		
		4 &      82.52 & {\bf 82.21} &      82.46 &      82.32 &      82.48 &     169.13 & {\bf 168.77} &     169.05 &     169.36 &     169.26 \\
		
		5 &      83.85 & {\bf 83.71} &      83.82 &      83.82 &      83.83 &     170.66 &     170.60 &     170.61 &     170.71 & {\bf 170.53} \\
		
		6 &      82.61 & {\bf 82.51} &      82.58 &      83.36 &      83.14 &     170.56 & {\bf 170.39} &     170.45 &     170.82 &     170.82 \\
		
		7 &      87.04 & {\bf 87.01} &      87.02 &      87.34 &      87.21 &     172.94 & {\bf 172.67} &     172.76 &     173.09 &     173.02 \\
		
		8 &      82.09 & {\bf 81.85} &      81.87 &      81.85 &      81.95 &     167.82 &     167.52 &     167.71 & {\bf 167.03} &     167.48 \\
		
		9 &      82.59 & {\bf 82.47} &      82.52 &      82.54 &      82.57 &     171.42 & {\bf 171.31} &     171.36 &     171.32 &     171.52 \\
		
		10 &      83.55 &      83.52 & {\bf 83.51} &      83.71 &      83.65 &     169.45 &     168.97 &     169.31 &     169.04 & {\bf 168.82} \\
		\bottomrule
	\end{tabular}

	\label{tab:results_mode}
\end{table*}
The expected values for electricity and natural gas prices are obtained by calculating a linear combination of the observations of the last five years weighted by the probability ($\widehat{x_{t}}=\sum_{i=1}^{5}\pi_{\omega_i}x_{t,i}$ where $x_{t}$ is the price for time period $t \in \mathcal{T}$ in year $i$). Due to the weekly time periods, the values are averaged per week. 

\subsection{Operational planning problem}
In the operational planning more recent information is available for the scenario generation, because we obtain new observations after each week. Furthermore, we are closer to actual delivery time than in the biomass contract selection problem. Consequently, we can use  time series analysis to better predict the uncertainties by updating the models in every week. 

There a different possibilities to obtain scenarios for the operational model. We implement and analyze five different types of scenario generation:
\paragraph{Using past data as predictions (P)} Data from previous years is used to built scenarios analog to the biomass contract selection scenarios. The scenarios consist of the data from the respective week(s) in previous years. 
\paragraph{Combining time series models and past data as predictions (F1)} In this method, we use time series models to predict the first week of the receding horizon and use data from previous years for the remaining weeks of the receding horizon.
The time series model uses the most recent observations to update the forecast for the following week. We use an ARMAX model \citep{madsen2007time} with weekly seasonality of prices and consumption using Fourier series in the form of exogenous parameters \citep{weron2007modeling}. We use past data for the remaining weeks because they are further into the future and the risk of inaccurate predictions is higher. 
To create scenarios from the time series model, we follow the scenario generation process described in \citep{conejo2010decision}. More specifically, we generate 2500 equiprobable scenarios using Monte Carlo simulation and cluster them using the k-medoid algorithm to obtain five  representative scenarios \citep{hastie2009unsupervised}. The forecasted scenarios for the first week have to be combined with data from previous years to get a scenario for the entire receding horizon. Therefore, we add the data from the most recent year to the scenario with the highest probability. 
\paragraph{Using time series models as predictions (F2)} This method is similar to \emph{F1}, because it also uses time series models for predictions and uses Monte Carlo simulation and clustering for generating scenarios. However, in this case we make predictions for the entire receding horizon and do not combine with past data.  The time series models, forecasts and scenarios are obtained following the same method as for \emph{F1}.  \\

All three above mentioned methods result in five scenarios for the operational planning problem. As two further possibilities for scenario generation, we use combinations of these methods. Namely, we combine the scenarios obtained from historical data (\emph{P}) with the two time-series-based methods (\emph{F1} and \emph{F2}) resulting in ten scenarios. Note that the probabilities are normalized to result in a sum of one again. These methods are denoted by \emph{P}+\emph{F1} and \emph{P}+\emph{F2}.

Note that the above mentioned scenario generation is used for electricity prices and heat demands. For the gas prices in case study A with the gas boiler, we also use an expected value in the operational model. This is due to the fact, that gas prices are daily prices and are not as volatile as, e.g., the electricity price, and we deem the expected value as accurate enough for this model.

\section{Analysis of  scenario generation methods}\label{sec:appendix_b}

In this section, we compare the different methods for scenario generation. The results show the performance of the scenario method configuration for  the operational  planning problem presented in Section \ref{sec:appendix_a} of this Appendix.

Table \ref{tab:results_mode} shows the results for each sample of both municipalities. The value shown is the minimum overall costs per scenario generation method, where the minimum is taken over the minimum objective value obtained for four different lengths of the receding horizon (one, two, three or four weeks). The analysis of different receding horizon lengths is described in Section 6.1 of the main article. Based on Table \ref{tab:results_mode}, the best of the implemented scenario generation methods is \emph{F1}, i.e., updating the scenarios every week by forecasting the next week of heat demand and using previous years for the remaining weeks of the receding horizon. Method \emph{F1} achieves the best result in 9 out of 11 samples for municipality A and in 8 out of 11 cases for municipality B. For the remaining 2 and 3 cases, respectively, no common favorable can be determined, as it differs per case. However, in all cases using the scenarios of method \emph{F1} is better than using expected values, as we show in Section 6.2 of the main article.

Based on these results, we conclude for our test cases that it is beneficial to update the scenarios every week instead of using previous years' data. However, using time series models for more than the first week often leads to worse results, which means that the scenarios are misleading the optimization. Therefore, using updated information just for the first week is a compromise and improves the results. For application in practice, this should be evaluated individually. Furthermore, our scenario generation methods can be easily replaced with already existing proved and tested forecasting methods of the operator.

\end{appendices}

\bibliography{Bibliography}

\end{document}